\documentclass[hidelinks,onefignum,onetabnum]{siamart251216}
\usepackage{xspace}
\usepackage{amsmath,amssymb,amsfonts,amsopn}
\usepackage{graphicx,tikz,algorithmic,xcolor}
\usetikzlibrary{arrows.meta}
\usepackage{array,booktabs,tabularx,multirow,mathtools,needspace,enumitem}

\newcommand{\hypre}{{\sl hypre}\xspace}
\colorlet{blue}{black}
\renewcommand{\footercopyright}{}
\hypersetup{
  pdftitle={Online Autotuning of Algebraic Multigrid via Contextual Bandits and Reinforcement Learning},
  pdfauthor={Xudong He, Ruiyu Wang, Yilang Ding, Rui Peng Li, and Yuanzhe Xi}
}

\headers
{Autotuning AMG}
{X. He, R. Wang, Y. Ding, R. Li, and Y. Xi}

\title{Online Autotuning of Algebraic Multigrid via Contextual Bandits
and Reinforcement Learning}

\author{%
Xudong He%
\thanks{%
Department of Mathematics, Emory University
(\email{xudong.he@emory.edu},
 \email{jonathan.wang@emory.edu},
 \email{lance.ding@emory.edu},
 \email{yxi26@emory.edu}).
The work of Y. Xi is supported by the National Science Foundation
under awards DMS-2038118 and DMS-2338904.
}%
\and
Ruiyu Wang\footnotemark[1]
\and
Yilang Ding\footnotemark[1]
\and
Rui Peng Li%
\thanks{%
Center for Applied Scientific Computing,
Lawrence Livermore National Laboratory,
P.O. Box 808, L-561, Livermore, CA 94551
(\email{li50@llnl.gov}).
This work was performed under the auspices of the
U.S. Department of Energy by Lawrence Livermore National Laboratory
under Contract DE-AC52-07NA27344.
LLNL-JRNL-XXXXXXX.
The United States Government retains, and the publisher, by accepting
the article for publication, acknowledges that the United States
Government retains a non-exclusive, paid-up, irrevocable, worldwide
license to publish or reproduce the published form of this manuscript,
or allow others to do so, for United States Government purposes.
}%
\and
Yuanzhe Xi\footnotemark[1]
}

\DeclareMathVersion{figureml}
\SetSymbolFont{operators}{figureml}{OT1}{mlmr}{m}{n}
\SetSymbolFont{letters}{figureml}{OML}{mlmm}{m}{it}
\SetSymbolFont{symbols}{figureml}{OMS}{mlmsy}{m}{n}
\SetSymbolFont{largesymbols}{figureml}{OMX}{mlmex}{m}{n}
\SetMathAlphabet{\mathbf}{figureml}{OT1}{mlmr}{bx}{n}
\SetMathAlphabet{\mathsf}{figureml}{OT1}{mlmss}{m}{n}
\SetMathAlphabet{\mathit}{figureml}{OT1}{mlmr}{m}{it}
\SetMathAlphabet{\mathtt}{figureml}{OT1}{mlmtt}{m}{n}
\newcommand{\MLModernFigureFont}{%
  \renewcommand{\rmdefault}{mlmr}%
  \renewcommand{\sfdefault}{mlmss}%
  \renewcommand{\ttdefault}{mlmtt}%
  \fontfamily{mlmr}\selectfont\mathversion{figureml}%
}

\newcommand{\ModuleFourSeedCoverage}{six seed settings for every family and grid}
\newcommand{\ModuleFourGrid}{60}
\newcommand{\ModuleFourSeed}{3}
\newcommand{\ModuleFourSeedCountWord}{six}
\newcommand{\ModuleFourDiffusionFullSavings}{36.1\%, 43.3\%, and 42.5\%}
\newcommand{\ModuleFourDiffusionFrozenSavings}{47.0\% to 49.8\%}
\newcommand{\ModuleFourAdvectionFullSavings}{19.8\%, 31.0\%, and 30.2\%}
\newcommand{\ModuleFourAdvectionFrozenSavings}{36.8\% to 40.5\%}
\newcommand{\ModuleFourSelectedDefaultCost}{293.4}
\newcommand{\ModuleFourSelectedJointCost}{160.4}
\newcommand{\ModuleFourSelectedFullReduction}{45.3}
\newcommand{\ModuleFourSelectedFrozenReduction}{50.6}
\newcommand{\ModuleFourExampleIDs}{1563, 1564, and 1565}
\newcommand{\ModuleFourExampleCycles}{11, 12, and 16}
\newcommand{\ModuleFourExampleGrid}{80}
\newcommand{\ModuleFourExampleSeed}{5}
\newcommand{\ModuleFourExampleInitialWeights}{2.05, 2.35,\text{ and }3.00}

\begin{document}

\maketitle

\begin{abstract}
Algebraic multigrid (AMG) methods are widely used for solving large sparse linear systems, but their efficiency depends strongly on interacting setup and solve parameters. Effective parameter choices typically depend on the problem. We develop a phase-aware online autotuning method for \hypre's BoomerAMG to reduce cumulative time-to-solution over sequences of related linear systems. A LinUCB-based contextual bandit selects mixed-type parameters for AMG hierarchy construction using a shared cost model over PDE descriptors and configuration encodings. A recursive LSTDQ controller uses residual and timing observations to estimate the AMG solve time within each attempt and adapt relaxation weights between multigrid cycles. Both learners update
from ongoing solves without offline pretraining. Across grids of sizes \(40^3\), \(60^3\), and \(80^3\), the joint autotuning method reduces total runtime, including recurring learner computation, by 36.1--43.3\% for diffusion problems and 19.8--31.0\% for diffusion--advection problems over sequences of 5000 linear systems, relative to the reference BoomerAMG default configuration. On a variable-flow diffusion--advection benchmark comprising 30,000 systems on a \(60^3\) grid, 2 systems fail to converge with joint autotuning, compared with 1677 for the reference configuration, under their respective attempt budgets.
\end{abstract}



\section{Introduction}
\label{sec:introduction}

Many scientific and engineering applications, including time-dependent partial differential equations, inverse problems, uncertainty quantification, and optimization loops, require the repeated solution of related sparse linear systems \cite{dambra2021amg,xu2020setup}. We represent such a sequence of $T$ linear systems as
\begin{equation}
    \mathbf{A}^{(t)}\mathbf{u}^{(t)}=\mathbf{f}^{(t)},
    \qquad t=1,2,\ldots,T,
    \label{eq:problem}
\end{equation}
where \(\mathbf{A}^{(t)}\in\mathbb{R}^{n_t\times n_t}\) is large and sparse, and \(t\) indexes the problem instance. The practical objective is to reduce cumulative time-to-solution across the sequence. Solving a sequence of related problems creates an opportunity for online autotuning: information collected from earlier solves can guide parameter choices on later instances.

We focus on algebraic multigrid (AMG) methods, which are scalable and can achieve near-linear complexity for many elliptic problems. These methods apply to matrices from unstructured discretizations or algebraically defined operators \cite{henson2002boomeramg,ruge1987algebraic}. Their practical efficiency depends strongly on setup and solve parameters, including the choices of coarsening, interpolation, and smoothing methods, and cycle structure, which affect convergence, operator complexity, memory use, communication, and runtime \cite{ruge1987algebraic,henson2002boomeramg,desterck2006reducing, baker2011multigrid,yang2010long}. \hypre's BoomerAMG is a scalable parallel AMG implementation. It exposes these choices through interacting continuous, integer-valued, and categorical parameters. Default configurations provide useful general-purpose baselines, but need not be the most efficient ones for a particular problem sequence.

Tuning AMG requires balancing hierarchy-construction cost against subsequent solve cost. Setup parameters determine the multilevel hierarchy, while solve-phase parameters control its application through relaxation and coarse-grid correction. Together, they determine the cost of each cycle and the number of cycles needed to reach tolerance. An AMG hierarchy that is inexpensive to construct or apply may require more cycles, so minimizing setup cost or iteration count alone need not minimize total solution time. Although AMG theory provides valuable insight into convergence and algorithm design \cite{henson2002boomeramg,maclachlan2014bounds,xu2017algebraic}, predicting these tradeoffs across a large, mixed-type parameter space remains challenging \cite{caldana2024deep,gahvari2011modeling}. This motivates learning from measured setup and solve costs.

Existing autotuners can exploit such measurements. GPTune shares information across tuning tasks using multi-task Gaussian-process (GP) models, while GPTuneBand additionally leverages evaluations at multiple fidelities \cite{liu2021gptune,zhu2022gptuneband}. These frameworks have been applied to multi-parameter tuning of \hypre. Online parameter selection for sequences of linear systems is studied in \emph{Learning to Relax} \cite{khodak2024learning}, which develops bandit and contextual-bandit methods for selecting SOR and SSOR relaxation parameters, together with learning-theoretic guarantees under the stated assumptions. Since SOR and SSOR can also be used as AMG smoothers, this work is closely related to our setting. However, it selects a single scalar relaxation weight for each system and keeps that weight fixed throughout the solve, rather than adapting a sequence of actions across multigrid cycles within a solve.

Neural networks have been used to learn AMG prolongation operators \cite{luz2020learning}, smoothers \cite{huang2023learning}, and sparse coarse-grid operators \cite{huang24}. Other approaches retain existing components and
learn their configuration. Caldana et al.\
\cite{caldana2024deep} use an offline-trained runtime surrogate
to select the AMG strength threshold. Antonietti et al.\
\cite{antonietti2026deep} extend this approach to jointly select the strength-of-connection (SoC) threshold and smoother type, accounting for interactions between hierarchy and relaxation choices. Both approaches select configurations using previously trained predictors. AMG parameters have also been adapted during solver execution. Mishra et al.\ \cite{mishra2024online} use online Bayesian optimization to update polynomial-multigrid cycle configurations
as a simulation progresses. Huergo et al.\
\cite{huergo2024multigrid} use proximal policy optimization (PPO) to adjust smoothing sweeps and a correction fraction between V-cycles in h/p-multigrid. Their controller responds to residual information and uses rewards based on residual reduction per unit time; the reported tests evaluate a previously trained policy.

Our contribution is a phase-aware online autotuning method that learns mixed-type setup configurations and cycle-level relaxation weights while solving a finite sequence of related linear systems. The method combines two components: a \emph{setup learner} for AMG setup parameters that control hierarchy construction and a \emph{solve learner} for AMG solve parameters, in particular the relaxation weights. They maintain separate models but are coupled through solver execution: the selected hierarchy affects solve performance, while the relaxation policy affects the cost and success of a setup choice. The distinction from prior work is the integration of these decisions with continued learning within the workload being solved.

The setup learner makes one configuration-level decision before each attempt. It observes feedback only for the selected configuration, including the cost of the subsequent solve, without modeling individual cycle transitions. We use a LinUCB-based contextual-bandit method \cite{li2010contextual,chu2011contextual} with a shared linear cost model over PDE descriptors and setup parameters. This representation allows observations from one configuration to inform predictions for related, untested configurations. The solve learner instead makes a sequence of decisions on the constructed hierarchy. The relaxation weight selected for each cycle affects the residual state on which subsequent cycles act, so maximizing immediate residual reduction does not necessarily minimize total solve time. We use recursive least-squares temporal-difference action-value estimation (LSTDQ) \cite{lagoudakis2003least} to predict the execution time of the remaining cycles within the current attempt. These predictions guide relaxation-weight selection with exploration, and the learner updates using cycle-level residual and timing feedback.

Both learners update from the first problem in the sequence without offline pretraining. Our main online comparisons include cumulative setup and solve time, parameter-selection and model-update costs, and every executed attempt, including failures and retries. The evaluation therefore asks whether solver savings outweigh learning costs within the available sequence, rather than measuring performance only after training. Experiments on structured anisotropic diffusion and diffusion--advection problems, including a spatially varying advection benchmark, evaluate runtime and convergence reliability under the prescribed retry and fallback protocols. Matched-hierarchy tests separately compare learned relaxation with grid-searched fixed weights and a prescribed period-two policy, distinguishing relaxation performance from changes in hierarchy selection. The results demonstrate runtime savings for the complete method on the tested sequences and improvements over the tested fixed-weight references.

The remainder of the paper is organized as follows. Section~\ref{sec:problem} defines the autotuning problem. Sections~\ref{sec:setup_bandit} and~\ref{sec:solve_rl} present the setup and solve learners, respectively.
Section~\ref{sec:experiments} reports the numerical results, and Section~\ref{sec:conclusions} concludes.

\section{Problem Setup and Formulation}
\label{sec:problem}

This section describes how the online tuning method processes each linear system \eqref{eq:problem}. We distinguish the AMG hierarchy-construction attempts from the AMG cycles within an attempt, define the information available to the learners, and specify the runtime objective used to evaluate the complete method.

\subsection{AMG parameters and execution}
\label{subsec:boomeramg_setup_solve}

BoomerAMG first constructs a multilevel hierarchy and then applies multigrid cycles using the hierarchy. We therefore partition its parameters as
\begin{equation}
    \Theta
    =
    \Theta_{\mathrm{setup}}
    \times
    \Theta_{\mathrm{solve}},
    \label{eq:theta_split}
\end{equation}
where setup parameters determine hierarchy construction and solve parameters determine how the cycle is performed.

For problem \(t\), let \(j\) index successive \emph{attempts}. Attempt \((t,j)\) first constructs a hierarchy using \(a_{t,j}\in\Theta_{\mathrm{setup}}\). If the construction succeeds, AMG cycles are performed using that hierarchy. The attempt fails either if hierarchy construction fails or if the following solve does not reach the prescribed tolerance within the cycle limit.
For methods with learned setup selection, a failed attempt may be followed by another learner-selected configuration while the attempt budget permits. If all learned attempts fail, methods with fallback make one final attempt using the reference setup and solve settings. This process stops at the first successful attempt.

Let \(J_{\max}\) denote the maximum number of attempts and \(J_t\le J_{\max}\) the number actually executed on problem \(t\). Method-specific attempt budgets are given in Section~\ref{subsec:experiment_setup}.
For a successful setup configuration \(a\), let \(\mathcal H_t(a)\) denote the resulting AMG hierarchy. Attempt \((t,j)\) uses \(\mathcal H_t(a_{t,j})\) throughout its solve. Changing solve phase parameters such as the relaxation weight between cycles does not rebuild the hierarchy; a new hierarchy is constructed only when a new setup attempt begins.
Let \(\mathbf u_k^{(t,j)}\) denote the iterate after \(k\) completed cycles of attempt \((t,j)\). Define the corresponding residual and relative residual norm reduction as
\begin{equation}
    \mathbf r_k^{(t,j)}
    =
    \mathbf f^{(t)}
    -
    \mathbf A^{(t)}\mathbf u_k^{(t,j)},
    \qquad
    \varrho_{t,j,k}
    =
    \frac{\|\mathbf r_k^{(t,j)}\|_2}
         {\|\mathbf r_0^{(t,j)}\|_2},
    \label{eq:residual}
\end{equation}
for a nonzero initial residual $\mathbf r_0^{(t,j)}$.
Let \(w_{t,j,k}\) denote the relaxation weight used in the cycle that advances \(\mathbf u_k^{(t,j)}\) to \(\mathbf u_{k+1}^{(t,j)}\). When solve-phase learning is enabled, the solve learner selects the weight using the information observed up to the $k$-th cycle.
Let \(K_{\max}\) be the maximum number of cycles per attempt and \(\varepsilon_{\rm tol}>0\) the relative residual tolerance. An attempt succeeds if there exists a $k$ such that \(\varrho_{t,j,k}<\varepsilon_{\rm tol}\) and $k < K_{\max}$;
otherwise, it is recorded as a failure.

Table~\ref{tab:boomeramg-params} lists representative BoomerAMG parameters and their reference settings. Our experiments tune the 8 setup controls in Table~\ref{tab:search_space} and the solve-phase relaxation weight; all remaining controls are fixed.

\begin{table}[t]
\centering
\caption{Representative BoomerAMG parameters and reference settings used in this study. Here \(r\), \(i\), and \(c\) denote real-valued, integer-valued, and categorical parameters, respectively.}
\label{tab:boomeramg-params}
\footnotesize
\setlength{\tabcolsep}{3.2pt}
\begin{tabular}{@{}p{0.2\linewidth}p{0.11\linewidth}
                  p{0.5\linewidth}p{0.10\linewidth}@{}}
\toprule
Parameter & Type & Role & Reference \\
\midrule
\multicolumn{4}{@{}l}{\textit{Setup parameters}} \\
\addlinespace[1pt]

\texttt{strong\_threshold}
& \(0\le r<1\)
& Strength-of-connection threshold.
& \(0.25\) \\

\texttt{max\_row\_sum}
& \(0<r<1\)
& Row-sum criterion related to diagonal dominance.
& \(0.9\) \\

\texttt{coarsen\_type}
& \(c\)
& Coarsening algorithm.
& \(10\) \\

\texttt{interp\_type}
& \(c\)
& Interpolation algorithm.
& \(6\) \\

\texttt{trunc\_factor}
& \(r\ge0\)
& Interpolation truncation factor.
& \(0.0\) \\

\texttt{P\_max\_elmts}
& \(i>0\)
& Maximum interpolation entries per row.
& \(4\) \\

\texttt{agg\_num\_levels}
& \(i\ge0\)
& Number of aggressive coarsening levels.
& \(0\) \\

\texttt{agg\_interp\_type}
& \(c\)
& Interpolation type with aggressive coarsening.
& \(4\) \\

\texttt{nongalerkin\_tol}
& \(r\ge0\)
& Non-Galerkin dropping tolerance.
& \(0.0\) \\

\midrule
\multicolumn{4}{@{}l}{\textit{Solve parameters}} \\
\addlinespace[1pt]

\texttt{relax\_type}
& \(c\)
& Relaxation method in the AMG cycles.
& \(18/18\) \\

\texttt{num\_sweeps}
& \(i>0\)
& Number of relaxation sweeps.
& \(1\) \\

\texttt{relax\_weight}
& \(r>0\)
& Weight used by the relaxation scheme.
& \(1.0\) \\

\texttt{cycle\_type}
& \(c\)
& Multigrid cycle type, such as a V-cycle or W-cycle.
& V \\
\bottomrule
\end{tabular}
\end{table}

\subsection{Online objective and runtime accounting}
\label{subsec:sequence_objective}

When processing problem~\(t\), the learners observe a context vector
$\chi_t\in\mathbb R^{d_\chi}$,
containing problem descriptors available before the AMG hierarchy construction. The same context is provided to both learners and remains fixed across all attempts on the problem. The specific descriptors used in the experiments are given in Section~\ref{subsec:experiment_setup}.

Let \(\Pi\) denote a complete online learning method, including its setup and solve decisions, retry and fallback rules, and learner updates. All decisions use only the information available at the time they are made. Let \(\tau_{t,j}^{\rm set}\) and \(\tau_{t,j}^{\rm sol}\) denote the setup and solve execution times of attempt \((t,j)\), excluding learner computation. If hierarchy construction fails, \(\tau_{t,j}^{\rm sol}=0\), while the time spent on construction remains included.

Let \(O_t^{\rm set}\) and \(O_t^{\rm sol}\) denote the total setup- and solve-learner computation recorded on problem \(t\), including action selection, target construction, and model updates. We then define
\begin{equation}
    \begin{aligned}
    C_t^{\rm solver}(\Pi)
    &=
    \sum_{j=1}^{J_t}
    \left(
        \tau_{t,j}^{\rm set}
        +
        \tau_{t,j}^{\rm sol}
    \right),\\
    C_t^{\rm total}(\Pi)
    &=
    C_t^{\rm solver}(\Pi)
    +
    O_t^{\rm set}
    +
    O_t^{\rm sol}.
    \end{aligned}
    \label{eq:timing_decomposition}
\end{equation}
Here \(C_t^{\rm solver}\) measures the setup and solve time, while \(C_t^{\rm total}\) additionally includes the time for learner computation, and \(O_t^{\rm set}=\sum_j O_{t,j}^{\rm set}\) and \(O_t^{\rm sol}=\sum_j O_{t,j}^{\rm sol}\) denote the setup- and solve-learner computation times, respectively.
For a sequence of \(T\) problems, the cumulative recorded runtime is
\begin{equation}
    \mathcal C_T^{\rm total}(\Pi)
    =
    \sum_{t=1}^{T} C_t^{\rm total}(\Pi).
    \label{eq:sequence_objective}
\end{equation}
The online objective is to reduce this cumulative cost, including learner computation from the beginning of the sequence.

Because runtime depends on the problem sequence, learner randomization, and timing variability, the experiments compare multiple independent runs. We also report the number of problems that fail to converge within the allowed attempt budgets, since runtime comparisons alone are insufficient when the compared solver configurations successfully solve different numbers of problems. Table~\ref{tab:notation} summarizes the main notation used throughout the paper.

\begin{table}[t]
\centering
\small
\caption{Main notation for problem indices, parameter choices, and runtime accounting.}
\label{tab:notation}
\begin{tabular}{@{}p{0.15\linewidth}>{\footnotesize}p{0.8\linewidth}@{}}
\toprule
Symbol & \small Meaning \\
\midrule

\(t,\ T\)
& Problem index and number of problems in the sequence. \\

\(j,\ J_t,\ J_{\max}\)
& Attempt index, number of attempts on problem \(t\),
and maximum attempt budget. \\

\(k,\ K_{\max}\)
& Number of completed cycles within an attempt and the
per-attempt cycle limit. \\

\(a_{t,j}\)
& Setup configuration selected for attempt \((t,j)\). \\

\(w_{t,j,k}\)
& Relaxation weight applied to cycle \(k+1\) of attempt \((t,j)\). \\

\(\varrho_{t,j,k}\)
& Relative residual after \(k\) completed cycles of
attempt \((t,j)\). \\

\(\chi_t\)
& Problem context available to both learners. \\

\midrule

\(\tau_{t,j}^{\rm set},\ \tau_{t,j}^{\rm sol}\)
& Setup and solve times for attempt \((t,j)\),
excluding learner computation. \\

\(O_t^{\rm set},\ O_t^{\rm sol}\)
& Setup- and solve-learner computation times totaled over problem \(t\). \\

\(C_t^{\rm solver},\ C_t^{\rm total}\)
& Solver execution time and total recorded runtime for
problem \(t\). \\

\(\mathcal C_T^{\rm total}\)
& Total recorded runtime over the \(T\)-problem sequence. \\

\bottomrule
\end{tabular}
\end{table}

\section{Setup-Phase Contextual Bandit}
\label{sec:setup_bandit}
The setup learner uses the problem context $\chi_t$ and observations from previously executed configurations to select an AMG-setup configuration. Since only the selected configuration is executed, the learner receives contextual-bandit feedback \cite{khodak2024learning,li2010contextual}. A context--action cost model, trained from previously executed attempts, ranks candidate configurations for the first setup attempt on each problem. If an attempt fails, its binary failure outcome is available immediately, whereas its final cost target may depend on subsequent attempts and is therefore not yet known. We consequently use a separate failure model, updated immediately after the failed attempt, to rank configurations for the next retry.

\subsection{Setup actions and feedback}
\label{subsec:setup_actions_loss}

Let \(\mathcal U\subseteq\Theta_{\mathrm{setup}}\) denote the finite setup search space.
Each configuration \(a=(a_1,\ldots,a_8)\in\mathcal U\) specifies the 8 control parameters that are listed in Table~\ref{tab:search_space}. The first five parameters, $a_1$--$a_5$, are numeric and discretized. The remaining three, $a_6$--$a_8$, are categorical, which select the algorithmic choices in BoomerAMG's setup. In particular, when aggressive coarsening is disabled, the aggressive-interpolation parameter is inactive and fixed at its reference value to avoid duplicate configurations.

\begin{table}[t]
\centering
\caption{Setup search space for \(a=(a_1,\ldots,a_8)\). Coarsening and interpolation values are algorithmic labels. Aggressive-interpolation truncation and row-limit controls remain fixed at zero.}
\label{tab:search_space}
\footnotesize
\setlength{\tabcolsep}{4pt}
\begin{tabular}{@{}
    >{\centering\arraybackslash}p{0.09\linewidth}
    >{\raggedright\arraybackslash}p{0.34\linewidth}
    >{\raggedright\arraybackslash}p{0.49\linewidth}@{}}
\toprule
Symbol & Parameter & Values \\
\midrule
\(a_1\)
& \texttt{strong\_threshold}
& \(0,0.05,\ldots,0.95\) \\
\(a_2\)
& \texttt{max\_row\_sum}
& \(10^{-6},0.05,0.10,\ldots,0.95\) \\
\(a_3\)
& \texttt{trunc\_factor}
& \(0,0.05,\ldots,0.95\) \\
\(a_4\)
& \texttt{P\_max\_elmts}
& \(2,4,6,8,12,16\) \\
\(a_5\)
& \texttt{agg\_num\_levels}
& \(0,1,2\) \\
\(a_6\)
& \texttt{coarsen\_type}
& \(0,2,3,6,8,10\) \\
\(a_7\)
& \texttt{interp\_type}
& \(0,2,3,4,6,8,17,18\) \\
\(a_8\)
& \texttt{agg\_interp\_type}
& \(4\), when \(a_5=0\);
  \(4,6,7\), when \(a_5>0\) \\
\bottomrule
\end{tabular}
\end{table}

The first setup decision of the online run uses the default configuration. For each subsequent setup decision, the learner constructs a candidate set \(\mathcal U_{t,j}\subseteq\mathcal U\) containing at most 1024 configurations. This candidate set includes the default, the most recently selected configuration, and previously evaluated configurations with low setup-learning loss. Additional candidates are generated by changing one active parameter of selected configurations, drawing samples from the discrete search space, and mapping scrambled Sobol points to the numeric grids. Sobol sequences provide low-discrepancy points that cover the parameter space evenly; scrambling randomizes these points while preserving their low-discrepancy structure \cite{sobol1967distribution,owen1995randomly}. Duplicates are removed, and unused precomputed candidates fill any remaining slots.

The candidate generation uses only the information available at the current attempt. For a retry, configurations that have failed on the current problem are removed. The learner scores all candidates in \(\mathcal U_{t,j}\) but constructs and evaluates only the selected configuration, \(a_{t,j}\). Each executed learned attempt provides two forms of feedback: a \emph{cost feedback} that measures the computational cost of the selected configuration and is used by the setup cost model, and a binary \emph{failure feedback} that is used by a separate failure model to guide retries. We define these two signals as follows.

\paragraph{Cost feedback}
Each learned attempt \((t,j)\) is assigned a cost \(B_{t,j}\) used to construct the setup-learning target. This training cost is distinct from the measured per-problem runtime \(C_t^{\rm total}\) and is defined by
\begin{equation}
    B_{t,j}
    =
        (\tau_{t,j}^{\rm set}
    +
        O_{t,j,\rm sel}^{\rm set})
    +  (\tau_{t,j}^{\rm sol}
        +O_{t,j}^{\rm sol})    
    + \sum_{i=1}^{J_{t-1}}O_{t-1,i,\rm upd}^{\rm set}
        ,
    \label{eq:assigned_attempt_cost}
\end{equation}
where the first two groups account for setup execution and candidate selection, and solve execution and controller computation, respectively. The time required to update the setup learner's models for the current attempt, including delayed updates using its finalized cost target, is not yet known. We therefore charge every learned attempt on problem \(t\) the total setup-update time from problem \(t-1\), with \(J_0=0\). Hence \(B_{t,j}\) is a training-cost convention rather than the elapsed runtime of the attempt. Reported runtime counts learner computation only when it is incurred.

The default fallback uses the reference setup and solve settings and therefore receives no setup-selection, solve-learner, or carried setup-update charge:
\begin{equation}
    B_{t,j}
    =
    \tau_{t,j}^{\rm set}
    +
    \tau_{t,j}^{\rm sol}.
\end{equation}
The cost feedback for a learned attempt is
\begin{equation}
    \ell_{t,j}^{\rm set}
    =
    \begin{cases}
        B_{t,j},
        & \text{if attempt }(t,j)\text{ succeeds},\\
        B_{t,j} + B_{t,j+1} + \cdots + B_{t,j+s},
        & \text{if attempt }(t,j)\text{ fails},
    \end{cases}
    \label{eq:setup_loss_walltime}
\end{equation}
where \(s\) is the smallest positive integer for which attempt \((t,j+s)\) successfully constructs a hierarchy, regardless of its solve convergence. If no subsequent hierarchy construction, including that of the default fallback, succeeds within the budget, we take \(j+s=J_t\), so the feedback includes the costs of all remaining failed attempts. Thus, a successful attempt receives its own assigned cost, while a failed attempt, whether due to construction failure or nonconvergence, receives the cumulative assigned cost through the attempt \((t,j+s)\). This hierarchy-dependent training feedback is designed to avoid favoring configurations that fail quickly during construction. All executed attempts still contribute to the measured total per-problem runtime \(C_t^{\rm total}\).

\paragraph{Failure feedback}
Each learned attempt also produces a binary outcome
\begin{equation}
    f_{t,j}
    =
    \begin{cases}
        0, & \text{if attempt }(t,j)\text{ succeeds}\\
        1, & \text{if attempt }(t,j)\text{ fails}
    \end{cases}
    \label{eq:setup_failure_label}
\end{equation}
which is used by the auxiliary failure model when selecting a retry. Unlike the cost feedback of a failed attempt that depends on the costs of subsequent attempts, \(f_{t,j}\) is known immediately when attempt \((t,j)\) terminates. The failure model can therefore incorporate this observation before the next retry is selected. The two feedback signals serve different purposes:
\begin{itemize}[leftmargin=2em]
\item \(\ell_{t,j}^{\rm set}\) trains the cost model used for initial setup selection;
\item \(f_{t,j}\) trains the failure model used to avoid unreliable configurations during retries.
\end{itemize}
Both models and their shared feature representation are defined in Section~\ref{subsec:setup_linucb}.

\subsection{Shared context--action models and training}
\label{subsec:setup_linucb}

The setup learner maintains two linear models that use the same features but different responses. The cost model predicts the setup-learning loss \(\ell_{t,j}^{\rm set}\), while the failure model uses the binary outcome \(f_{t,j}\).
We first encode a setup configuration \(a\). The five numeric
controls are scaled as
\begin{equation}
    z_r
    =
    {(a_r-\bar a_r)}/{s_r},
    \qquad r=1,\ldots,5,
\end{equation}
where
$\bar a=(0.25,0.90,0,4,0)$, and
$s=(0.25,0.10,0.20,4,1)$.
We encode the categorical parameters using one-hot encoding, with the default category chosen as the reference. Each retained indicator equals \(1\) when its category is selected and \(0\) otherwise, so the default category is represented by an all-zero vector. For \(a_6\), \(a_7\), and \(a_8\), the default values are \(10\), \(6\), and \(4\), respectively. These parameters have \(6\), \(8\), and \(3\) categories, including the default, giving \(5\), \(7\), and \(2\) indicators, respectively. Concatenating the three indicator vectors gives \(c(a)\in\{0,1\}^{14}\).
With \(z=(z_1,\ldots,z_5)^\top\),
we encode the setup configuration by
\begin{equation}
    g(a)
    =
    \bigl(
        z^\top,\,
        (z_r^2)_{r=1}^{5},\,
        (z_rz_s)_{1\le r<s\le5},\,
        c(a)^\top,\,
        (c(a)\otimes z)^\top
    \bigr)^\top
    \in\mathbb R^{104},
\end{equation}
which includes the scaled numeric controls, their squares and pairwise products, the categorical indicators, and numeric--categorical interactions.

We then combine this setup encoding with the problem context. 
Let \(\chi \in\mathbb R^{d_\chi}\) denote a generic problem context vector,
where
\(\chi=(1,\chi_2,\ldots,\chi_{d_\chi})^\top\)
with \(\chi(1)=1\) providing an intercept. 
Define the feature vector used by both models as
\begin{equation}
    \phi^{\rm set}(\chi,a)
    =
    \bigl(
        \chi^\top,\,
        g(a)^\top,\,
        \chi_2g(a)^\top,\,
        \ldots,\,
        \chi_{d_\chi}g(a)^\top
    \bigr)^\top
    \in\mathbb R^{d_{\rm set}},
    \quad
    d_{\rm set}=105d_\chi.
\end{equation}
The interaction terms 
$\chi_ig(a)^\top$
allow the effect of a setup configuration to depend on the problem context.
Each executed learned attempt or successful default fallback provides one observation for both models. Index these observations by \(i=1,\ldots,n\), with \((t_i,j_i)\) denoting the corresponding problem and attempt. 
The two models use the same input \(\phi^{\rm set}(\chi_{t_i},a_{t_i,j_i})\), but different responses
that are
$y_i^{\rm set}
    =
    \ell_{t_i,j_i}^{\rm set}$ in \eqref{eq:setup_loss_walltime},
and
$f_i
    =
    f_{t_i,j_i}$ in \eqref{eq:setup_failure_label}, respectively.
\paragraph{Cost model}
The cost model is
\begin{equation}
    \widehat\mu_n^{\rm set}(\chi,a)
    =
    \phi^{\rm set}(\chi,a)^\top
    \widehat\theta_n^{\rm set}.
    \label{eq:shared_lin_model}
\end{equation}
With \(\lambda_{\rm set}=1\), define
\begin{equation}
    \begin{aligned}
    \mathsf V_n^{\rm set}
    &=
    \lambda_{\rm set}I
    +
    \sum_{i=1}^{n}
    \phi^{\rm set}(\chi_{t_i},a_{t_i,j_i})
    \phi^{\rm set}(\chi_{t_i},a_{t_i,j_i})^\top,\\
    \mathbf b_n^{\rm set}
    &=
    \sum_{i=1}^{n}
    y_i^{\rm set}
    \phi^{\rm set}(\chi_{t_i},a_{t_i,j_i}),\\
    \widehat\theta_n^{\rm set}
    &=
    (\mathsf V_n^{\rm set})^{-1}
    \mathbf b_n^{\rm set}.
    \end{aligned}
    \label{eq:shared_linucb_stats}
\end{equation}

\paragraph{Failure model}
The failure model uses the same feature matrix \(\mathsf V_n^{\rm set}\), but replaces the cost response \(y_i^{\rm set}\) by the failure response \(f_i\):
\begin{equation}
    \mathbf b_n^{\rm fail}
    =
    \sum_{i=1}^{n}
    f_i\,
    \phi^{\rm set}(\chi_{t_i},a_{t_i,j_i}),
    \qquad
    \widehat\theta_n^{\rm fail}
    =
    (\mathsf V_n^{\rm set})^{-1}
    \mathbf b_n^{\rm fail}.
\end{equation}
Its prediction for candidate \(a\) is
\begin{equation}
    \widehat f_n(\chi,a)
    =
    \phi^{\rm set}(\chi,a)^\top
    \widehat\theta_n^{\rm fail}.
\end{equation}
Because this is least-squares regression on binary responses, \(\widehat f_n\) is a ranking score rather than a calibrated failure probability and need not lie in \([0,1]\).

The two models are updated online as learned attempts are executed. Failure feedback is available immediately, whereas the final cost response may depend on subsequent attempts; the corresponding update procedure is described in Section~\ref{subsec:setup_selection}.

\subsection{LinUCB selection, retries, and online feedback}
\label{subsec:setup_selection}

Let \(n\) denote the number of setup-learning observations incorporated before the current selection. For context \(\chi\) and candidate \(a\), define the LinUCB width
\begin{equation}
    \sigma_n^{\rm set}(\chi,a)
    =
    \sqrt{
        \phi^{\rm set}(\chi,a)^\top
        (\mathsf V_n^{\rm set})^{-1}
        \phi^{\rm set}(\chi,a)
    }.
    \label{eq:setup_prediction}
\end{equation}
The first setup decision of the online run uses BoomerAMG's default. For the first attempt of each subsequent problem, candidates are ranked using the cost model:
\begin{equation}
    \mathrm{LCB}_n^{\rm set}(\chi,a)
    =
    \widehat\mu_n^{\rm set}(\chi,a)
    -
    \alpha_n\sigma_n^{\rm set}(\chi,a),
    \qquad
    \alpha_n=\frac{1}{\sqrt{n+1}},
\end{equation}
and the candidate with the smallest score is selected. Because cost is minimized, the negative uncertainty term temporarily favors less-observed feature directions. This effect decreases as \(\alpha_n\) decays with \(n\).

If a learned attempt fails and another learned attempt is allowed, its outcome is incorporated before the retry is selected. Configurations that have already failed on the current problem are removed, and the remaining candidates are ranked using the failure model:
\begin{equation}
    \mathrm{UCB}_n^{\rm fail}(\chi,a)
    =
    \widehat f_n(\chi,a)
    +
    \beta_{\rm fail}\sigma_n^{\rm set}(\chi,a),
    \qquad
    \beta_{\rm fail}=2.
    \label{eq:setup_failure_score}
\end{equation}
The candidate with the smallest score is selected. A smaller \(\widehat f_n\) indicates less association with previous failures, while the positive uncertainty term penalizes poorly observed feature directions. Thus, uncertainty encourages exploration for the first attempt but discourages uncertain choices after a failure.
Ties in \(\mathrm{UCB}_n^{\rm fail}\) are broken using the current cost LCB, followed by random selection if necessary. If the learned-attempt budget is exhausted, a final attempt uses the default setup and solve settings.
After each learned attempt, its feature vector is incorporated once into the shared matrix \(\mathsf V_n^{\rm set}\), and the failure response \(f_i\) is immediately incorporated into the failure model. A retry on the same problem therefore uses the newly observed failure information.

The cost model is also updated at this point, but the final response may not yet be available. If an attempt fails, \(\ell_{t,j}^{\rm set}\) depends on the costs of subsequent attempts through \eqref{eq:setup_loss_walltime}. The implementation therefore first uses the current cost prediction as a provisional response. Once the required subsequent cost is known, the provisional response is replaced by the final value \(y_i^{\rm set}=\ell_{t_i,j_i}^{\rm set}\). This correction changes the cost response but does not add a second observation to \(\mathsf V_n^{\rm set}\).

All executed learned attempts are retained as setup-learning observations. A default fallback enters the cost target of the preceding failed learned attempt. If it succeeds, it also provides a setup-learning observation for the default configuration, with cost target \(B_{t,j}\) and failure label zero. An unsuccessful default fallback does not add a separate setup-learning observation.

\section{Solve-Phase Reinforcement Learning}
\label{sec:solve_rl}

The setup learner selects the parameters used to construct the AMG hierarchy, whereas the solve learner selects a relaxation weight before each AMG V-cycle. The hierarchy remains fixed throughout an attempt, but the relaxation weight chosen at one cycle affects the residual state, and therefore the effectiveness and cost of subsequent cycles. Thus, the solve decisions are sequentially coupled. We therefore train the solve learner to minimize the expected execution time of the remaining cycles in the current attempt, rather than optimizing the residual reduction or execution time of the next cycle alone.

Sequential multigrid control has previously been studied using offline-trained PPO policies \cite{huergo2024multigrid}. In contrast, our solve learner updates online throughout the problem sequence and conditions its decisions on the problem context, selected setup configuration, and observed solve trajectory. We use an online variant of least-squares temporal-difference action-value estimation (LSTDQ) \cite{lagoudakis2003least}, augmented with eligibility traces following the LSTD($\lambda$) construction \cite{boyan2002technical}.

\subsection{Cycle-level control and remaining-time objective}
\label{subsec:solve_control_problem}

Fix problem \(t\) and attempt \(j\), with setup configuration \(a_{t,j}\) and hierarchy \(\mathcal H_t(a_{t,j})\). The hierarchy remains fixed throughout the attempt. During the solve, the learner repeatedly observes the current solver state, selects a relaxation weight, executes one V-cycle, and observes the resulting cycle cost and new state.

Let \(k\) denote the number of completed cycles and \(\mathbf u_k^{(t,j)}\) the current iterate. Before cycle \(k+1\), the learner selects
\begin{equation}
    w_{t,j,k}\in\mathcal W,
    \qquad
    \mathcal W=\{1,1.05,\ldots,3\}.
\end{equation}
All other solve settings remain fixed: one \(\ell_1\)-Jacobi sweep in each direction and a direct coarsest-level solve. The solve terminates at convergence, the cycle limit \(K_{\max}\), or solver failure. Let \(K_{t,j}\le K_{\max}\) denote the number of executed cycles.
The cost observed after applying \(w_{t,j,k}\) is
\begin{equation}
    c_{t,j,k}^{\rm sol}
    =
    \text{solver execution time of cycle }k+1.
    \label{eq:solve_cost}
\end{equation}
This includes solver-side parameter updates, the AMG cycle, and residual and stopping checks, but excludes solve-learner computation.
Let \(s_{t,j,k}\) denote the solver state immediately before cycle \(k+1\). It summarizes the fixed problem and hierarchy together with the solve history available at that point. To use this state in the LSTDQ model, we represent it by the feature vector \(h(s_{t,j,k})\) defined in \eqref{eq:solve_state} in Section~\ref{subsec:solve_training}.
Within an attempt, we suppress \((t,j)\), writing \(s_k\), \(w_k\), \(c_k^{\rm sol}\), and \(K\). Thus, each cycle produces the transition
\[
    (s_k,w_k)\longrightarrow(c_k^{\rm sol},s_{k+1}).
\]
The solve learner seeks weights that reduce the time required to complete the remaining cycles. To define this learning target, consider a specified, possibly randomized solve policy \(\pi^{\rm sol}\). Its action value is
\begin{equation}
    Q^{\pi^{\rm sol}}(s,w)
    =
    \mathbb E_{\pi^{\rm sol}}\!\left[
        \sum_{r=k}^{K-1}c_r^{\rm sol}
        \,\middle|\,s_k=s,\ w_k=w
    \right].
    \label{eq:solve_q_target}
\end{equation}
Hence \(Q^{\pi^{\rm sol}}(s,w)\) is the expected remaining cycle execution time after choosing weight \(w\) in state \(s\) and following \(\pi^{\rm sol}\) thereafter. Costs are undiscounted (\(\gamma=1\)), and the value after termination is zero. The target in \eqref{eq:solve_q_target} measures remaining cycle execution time until the current attempt terminates. It excludes learner computation, subsequent recovery attempts, and any terminal failure penalty.

\subsection{LSTDQ model and online training}
\label{subsec:solve_training}

The solve learner approximates the remaining-time target in \eqref{eq:solve_q_target}. For a specified solve policy,
\begin{equation}
    Q^{\pi^{\rm sol}}(s_k,w_k)
    =
    \mathbb E_{\pi^{\rm sol}}\!\left[
        c_k^{\rm sol}
        +
        Q^{\pi^{\rm sol}}(s_{k+1},w_{k+1})
        \,\middle|\,s_k,w_k
    \right],
\end{equation}
with zero successor value at termination. LSTDQ uses this one-step relation to learn a remaining-time predictor from observed cycle transitions.

\paragraph{State representation}
The solve learner represents the current state using the problem context, selected hierarchy, cycle index, and recent solve history. Within attempt \((t,j)\), write \(a=a_{t,j}\), \(\varrho_k=\varrho_{t,j,k}\), and \(\chi_t=(1,\xi^\top)^\top\), with \(\xi\in\mathbb R^{d_\chi-1}\). For \(k\ge1\), define
\begin{equation}
\begin{alignedat}{2}
    g_k
    &=
    \frac{\log(\varrho_k/\varepsilon_{\rm tol})}
         {\log(1/\varepsilon_{\rm tol})},
    &\qquad
    d_k
    &=
    \frac{\log(\varrho_k/\varrho_{k-1})}
         {\log(1/\varepsilon_{\rm tol})},\\
    q_k
    &=
    \log\!\left(1+\frac{c_{k-1}^{\rm sol}}{10^{-3}\,\mathrm{s}}\right),
    &\qquad
    \bar w_k
    &=
    w_{k-1}-2 .
\end{alignedat}
\end{equation}
These features measure the residual gap to tolerance, recent residual change, logarithmically scaled absolute cycle time, and previous relaxation weight. We initialize $(g_0,d_0,q_0,\bar w_0)=(1,0,0,-1)$.

To identify the AMG hierarchy used by the solve learner, encode the selected setup configuration by
\begin{equation}
    \zeta(a)
    =
    (\widetilde z^\top,c(a)^\top)^\top
    \in\mathbb R^{19},
    \qquad\widetilde z_r=\frac{2(a_r-a_r^{\min})}{a_r^{\max}-a_r^{\min}}-1,
\end{equation}
using the five numeric controls normalized to \([-1,1]\) over their allowed ranges and 14 categorical indicators from Section~\ref{subsec:setup_linucb}. This descriptor is fixed within an attempt. The state feature vector is
\begin{equation}
    h(s_k)
    =
    \bigl(
        1,\ k/K_{\max},\ g_k,\ g_k^2,\ d_k,\ q_k,
        \bar w_k,\ \xi^\top,\ \zeta(a)^\top
    \bigr)^\top
    \in\mathbb R^{25+d_\chi}.
    \label{eq:solve_state}
\end{equation}

\paragraph{Weight representation and predictor}
The action set contains 41 relaxation weights. To reduce the dimension of the LSTDQ model and the cost of evaluating all actions, we represent \(w\) using nine localized basis functions centered at \(\kappa_r=1+0.25(r-1)\), \(r=1,\ldots,9\):
\begin{equation}
    \begin{aligned}
    \widetilde b_r(w)
    &=
    \exp\!\left(
        -\frac{(w-\kappa_r)^2}{2(0.2)^2}
    \right)
    \mathbf 1\{|w-\kappa_r|\le0.5\},\\
    b_r(w)
    &=
    \frac{\widetilde b_r(w)}
         {\sum_{v=1}^{9}\widetilde b_v(w)},
    \qquad
    b(w)=(b_1(w),\ldots,b_9(w))^\top .
    \end{aligned}
    \label{eq:action_basis}
\end{equation}
Nearby weights therefore have similar representations.

Combining the action and state representations gives
\begin{equation}
    \phi^{\rm sol}(s,w)
    =
    b(w)\otimes h(s)
    \in\mathbb R^{d_{\rm sol}},
    \qquad
    d_{\rm sol}=9(25+d_\chi).
    \label{eq:solve_features}
\end{equation}
The remaining-time predictor is
\begin{equation}
    \widehat Q_n(s,w)
    =
    (\widehat\theta_n^{\rm sol})^\top
    \phi^{\rm sol}(s,w),
    \label{eq:q_features}
\end{equation}
where \(\widehat\theta_n^{\rm sol}\) is estimated by LSTDQ below. Writing \(\widehat\theta_n^{\rm sol}\) as nine blocks \(\widehat\theta_{n,r}^{\rm sol}\) gives
\[
\widehat Q_n(s,w)=\sum_{r=1}^9 b_r(w)\bigl(\widehat\theta_{n,r}^{\rm sol}\bigr)^\top h(s).
\]
Thus, for a fixed state, only the nine quantities \((\widehat\theta_{n,r}^{\rm sol})^\top h(s)\) need to be computed. They are then reused with the precomputed \(b(w)\) to evaluate all 41 candidate weights, rather than evaluating 41 separate state-dependent predictors.

\paragraph{LSTDQ update}
For transition \(i\), let
\begin{equation}
    \phi_i^{\rm sol}
    =
    \phi^{\rm sol}(s_i,w_i),
    \qquad
    \phi_i^{\rm sol\prime}
    =
    \phi^{\rm sol}(s_{i+1},w_{i+1}),
\end{equation}
where \(w_{i+1}\) is the next weight actually executed. At termination, \(\phi_i^{\rm sol\prime}=0\).
To propagate each observed cycle cost to recent state--weight pairs, LSTDQ maintains an eligibility-trace vector \(\mathbf z_i\):
\begin{equation}
    \mathbf z_i
    =
    \lambda_{\rm tr}\mathbf z_{i-1}
    +
    \phi_i^{\rm sol},
    \label{eq:eligibility_trace}
\end{equation}
where \(\lambda_{\rm tr}=0.8\) controls the geometric decay of contributions from preceding cycles. Thus, \(\mathbf z_i\) combines the current feature with decaying contributions from earlier transitions in the same solve attempt.

With regularization parameter \(\lambda_{\rm sol}=1\), the LSTDQ system is
\begin{align}
    \mathsf A_n^{\rm sol}
    &=
    \lambda_{\rm sol}I
    +
    \sum_{i=1}^{n}\mathbf z_i
    (\phi_i^{\rm sol}-\phi_i^{\rm sol\prime})^\top,
    \label{eq:lstdq}\\
    \mathbf b_n^{\rm sol}
    &=
    \sum_{i=1}^{n}\mathbf z_i c_i^{\rm sol},
    \nonumber\\
    \widehat\theta_n^{\rm sol}
    &=
    (\mathsf A_n^{\rm sol})^{-1}
    \mathbf b_n^{\rm sol}.
    \nonumber
\end{align}
The eligibility trace is reset between solve attempts: for the first transition of each attempt, the preceding trace is set to zero.

\subsection{Weight selection and exploration}
\label{subsec:solve_selection}

At each cycle, the solve learner uses the current LSTDQ predictor \(\widehat Q_n(s,w)\) together with a mismatch-based exploration score \(\sigma_{n,m}^{\rm sol}(s,w)\). With a small probability it instead selects a weight uniformly at random. We define these two exploration mechanisms below.

\paragraph{Mismatch-based exploration}
Let \(m\) denote the number of solve attempts retained in the LSTDQ model. For attempt \(e\), let \(\mathcal I_e\) denote its cycle transitions and define
\begin{equation}
    \Delta\mathsf A_e^{\rm sol}
    =
    \sum_{i\in\mathcal I_e}
    \mathbf z_i
    (\phi_i^{\rm sol}-\phi_i^{\rm sol\prime})^\top,
    \qquad
    \Delta\mathbf b_e^{\rm sol}
    =
    \sum_{i\in\mathcal I_e}
    \mathbf z_i c_i^{\rm sol}.
\end{equation}
Let \(\widehat\theta_e^{\rm end}\) be the fitted LSTDQ coefficient vector at the end of attempt \(e\), and define its mismatch by
\begin{equation}
    \mathbf g_e
    =
    \Delta\mathbf b_e^{\rm sol}
    -
    \Delta\mathsf A_e^{\rm sol}
    \widehat\theta_e^{\rm end}.
\end{equation}
Accumulating these mismatches gives
\begin{equation}
    \begin{aligned}
    \mathsf\Omega_m^{\rm sol}
    &=
    (10^{-3})^2\lambda_{\rm sol}I
    +
    \sum_{e=1}^{m}\mathbf g_e\mathbf g_e^\top,\\
    \sigma_{n,m}^{\rm sol}(s,w)^2
    &=
    \phi^{\rm sol}(s,w)^\top
    (\mathsf A_n^{\rm sol})^{-1}
    \mathsf\Omega_m^{\rm sol}
    (\mathsf A_n^{\rm sol})^{-\top}
    \phi^{\rm sol}(s,w).
    \end{aligned}
    \label{eq:cluster_score}
\end{equation}
A larger \(\sigma_{n,m}^{\rm sol}(s,w)\) indicates a state--weight direction associated with larger past LSTDQ mismatch. It is used as an exploration score rather than as a calibrated standard error.

The model-based weight is selected by
\begin{equation}
    w_k
    \in
    \operatorname*{arg\,min}_{w\in\mathcal W}
    \max\!\left\{
        \widehat Q_n(s_k,w)
        -2\sigma_{n,m}^{\rm sol}(s_k,w),\,0
    \right\}.
    \label{eq:solve_weight_choice}
\end{equation}
Thus, a small predicted remaining time is preferred, while a large mismatch score encourages exploration.

\paragraph{Random exploration}
To provide model-independent exploration, we also select a weight uniformly from \(\mathcal W\) with probability
\begin{equation}
    \epsilon_n^{\rm rand}
    =
    0.03+0.27\exp(-n/20000),
\end{equation}
where \(n\) is the total number of solve-cycle transitions observed so far. The first solve action uses \(w=1\); thereafter, \eqref{eq:solve_weight_choice} is used with probability \(1-\epsilon_n^{\rm rand}\). Hence, random exploration is strongest early in the online run and decreases as more transitions are observed.

\begin{figure}[t]
\centering
\begingroup\MLModernFigureFont
\resizebox{\linewidth}{!}{%
\begin{tikzpicture}[
  font=\footnotesize,
  box/.style={draw=black!55, rounded corners=2pt, line width=.55pt,
    align=center, text width=5.15cm, inner sep=6pt, fill=white},
  setup/.style={box, draw=blue!55!black, fill=blue!4},
  solve/.style={box, draw=teal!65!black, fill=teal!5},
  flow/.style={-{Stealth[length=2mm]}, line width=.65pt, draw=black!75},
  note/.style={font=\scriptsize, align=center, fill=white, inner sep=1.5pt}
]

\node[box, text width=11.8cm] (input) at (0,0)
  {\textbf{Problem sequence:} $\mathbf A^{(t)}\mathbf u^{(t)}=\mathbf f^{(t)}$,\quad $t=1,\ldots,T$
   \\Observe context $\chi_t$; start a new problem with attempt $j=1$.};

\node[font=\small\bfseries, text=blue!55!black] at (-3.45,-1.15)
  {SETUP: contextual bandit};
\node[font=\small\bfseries, text=teal!65!black] at (3.45,-1.15)
  {SOLVE: recursive LSTDQ};

\node[setup, minimum height=2.4cm] (select) at (-3.45,-2.55)
  {\textbf{Select setup $a_{t,j}$}\\[3pt]
   First attempt: $\arg\min_a(\widehat\mu^{\rm set}-\alpha_n\sigma^{\rm set})$\\[2pt]
   Retry: $\arg\min_a(\widehat f+2\sigma^{\rm set})$\\[2pt]
   Exclude setups that have already failed on problem $t$.};

\node[setup, minimum height=1.1cm] (build) at (-3.45,-4.55)
  {\textbf{Construct hierarchy $\mathcal H_t(a_{t,j})$}\\
   Measure setup execution time $\tau^{\rm set}_{t,j}$.};

\node[solve, minimum height=2.4cm] (weight) at (3.45,-2.55)
  {\textbf{Select relaxation weight $w_k\in\mathcal W$}\\[3pt]
   Explore uniformly, or minimize\\[2pt]
   $\max\{\widehat Q(s_k,w)-2\sigma^{\rm sol}(s_k,w),0\}$\\[2pt]
   $\widehat Q$: remaining cycle time\\in the current attempt.};

\node[solve, minimum height=1.1cm] (cycle) at (3.45,-4.55)
  {\textbf{Execute one V-cycle using $w_k$}\\
   Keep the hierarchy $\mathcal H_t(a_{t,j})$ fixed.};

\node[solve, minimum height=3.0cm] (observe) at (3.45,-6.95)
  {\textbf{Residual and timing feedback}\\[3pt]
   Update recursive LSTDQ.\\[2pt]
   Stop at convergence,\\cycle cap, or solver failure.};

\node[setup, minimum height=3.0cm] (feedback) at (-3.45,-6.95)
  {\textbf{Update setup models}\\[3pt]
   $\ell^{\rm set}_{t,j}=B_{t,j}+f_{t,j}\sum_{i=j+1}^{j+s} B_{t,i}$\\[2pt]
   $f_{t,j}=\mathbf{1}\{\text{attempt }(t,j)\text{ fails}\}$\\[2pt]
   $\widehat f_n(\chi,a)=\phi^{\rm set}(\chi,a)^\top
     \widehat\theta_n^{\rm fail}$\\[2pt]
   Failure feedback now; finalize cost later.};

\node[box, draw=orange!65!black, fill=orange!6,
      minimum height=2.4cm] (recover) at (-3.45,-10.35)
  {\textbf{Retry or BoomerAMG's\\default fallback}\\[3pt]
   Up to $J_{\max}-1$ learned attempts \\
   then one fixed Default setup and solve.\\
   Fallback cost updates\\the prior setup targets.};

\node[box, minimum height=2.4cm] (finish) at (3.45,-10.35)
  {\textbf{Finish problem $t$}\\[3pt]
   Solution or final failure.\\
   Record runtime; retain both models.};

\coordinate (entry) at (-5.95,0);
\draw[flow] (input.south -| entry) -- (select.north -| entry);
\draw[flow] (select) -- (build);
\draw[flow] (build.east) -- (0,-4.55) -- (0,-2.55)
  -- node[note, above=3pt, pos=0] {built\\$k=0$} (weight.west);
\draw[flow] (weight) -- (cycle);
\draw[flow] (cycle) -- (observe);
\draw[flow] (observe.east) -- (6.6,-6.95) -- (6.6,-2.55)
  -- (weight.east);
\node[note, rotate=90] at (6.82,-4.75)
  {continue: update $s_k$, choose next $w_k$};
\draw[flow] (observe.west) -- node[note, above=3pt] {attempt\\ends} (feedback.east);
\draw[flow] (build.west) -- (-6.6,-4.55) -- (-6.6,-6.95)
  -- (feedback.west);
\node[note, rotate=90] at (-6.82,-5.75) {setup failure};

\coordinate (outcome) at (-3.45,-8.75);
\draw[flow] (feedback.south) -- (outcome)
  -- node[note, right] {failure} (recover.north);
\draw[flow] (outcome) -- node[note, above=3pt] {success} (3.45,-8.75)
  -- (finish.north);
\draw[flow] (recover.east) -- node[note, above=3pt] {fallback\\ends} (finish.west);
\draw[flow] (recover.west) -- (-7.25,-10.35) -- (-7.25,-2.55)
  -- (select.west);
\node[note, rotate=90] at (-7.48,-5.5)
  {retry if $j<J_{\max}-1$: $j\leftarrow j+1$};
\draw[flow] (finish.east) -- (7.25,-10.35) -- (7.25,0) -- (input.east);
\node[note, rotate=90] at (7.48,-5.1)
  {next problem, if $t<T$: $t\leftarrow t+1$};

\end{tikzpicture}%
}
\endgroup
\caption{Online LinUCB--LSTDQ autotuning of BoomerAMG. The contextual bandit selects the AMG setup between solve attempts, while recursive LSTDQ selects the relaxation weight between V-cycles within a fixed hierarchy. Failed learned configurations trigger retries and, after the retry budget is exhausted, a fixed Default fallback. Both learners retain information across the problem sequence and optimize the cumulative recorded runtime.}
\label{fig:online-boomeramg-tuning}
\end{figure}
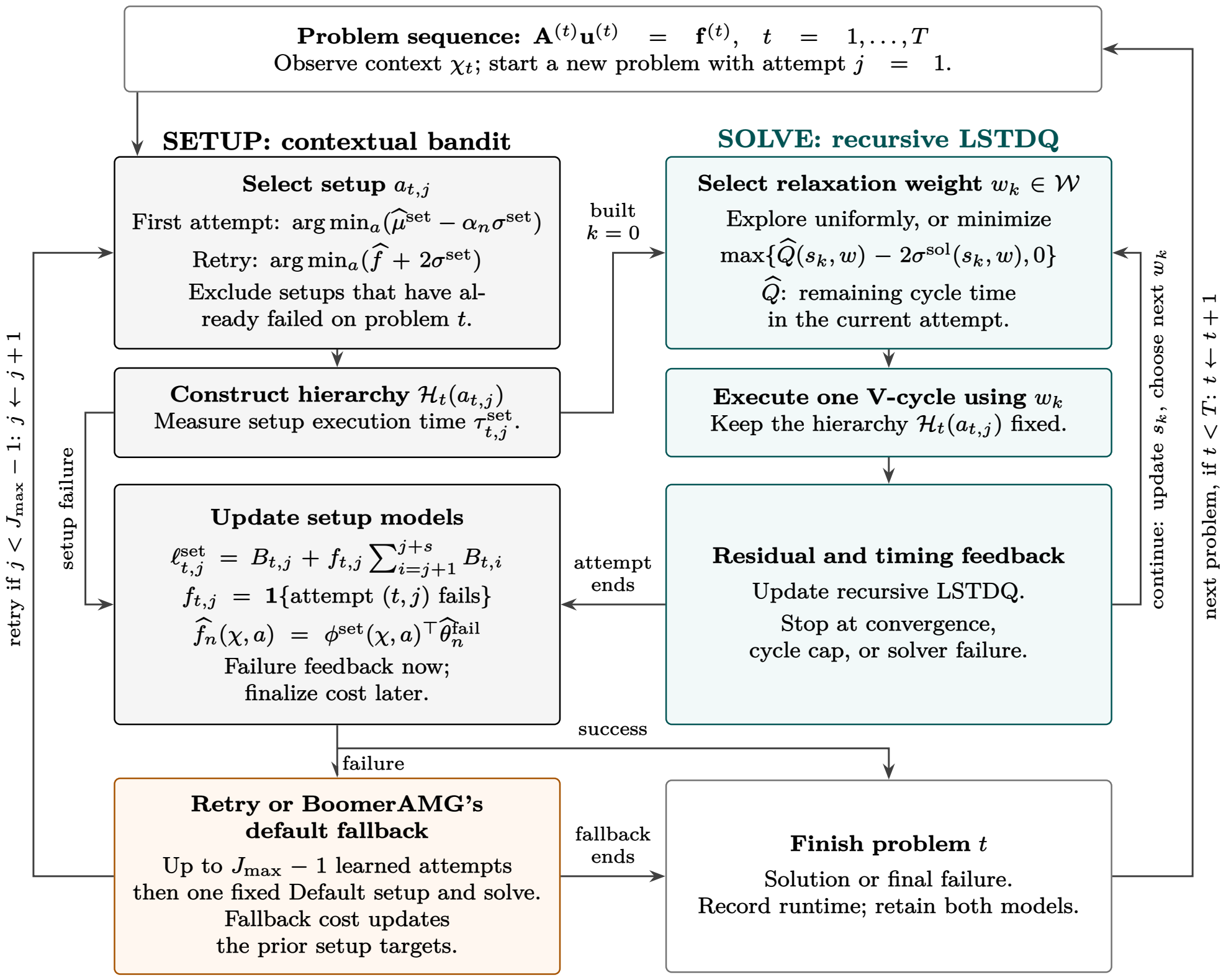

Figure~\ref{fig:online-boomeramg-tuning} summarizes the complete LinUCB--LSTDQ procedure. The setup learner selects a BoomerAMG configuration between attempts, while the solve learner selects relaxation weights between V-cycles within the resulting fixed hierarchy. Failed learned configurations may trigger another setup attempt, followed by the fixed default configuration after at most \(J_{\max}-1\) learned attempts. Both learners retain information across problems.

\section{Numerical Experiments}
\label{sec:experiments}

In this section, we evaluate the online autotuner on sequences of structured anisotropic diffusion and diffusion--advection systems. The experiments address four questions:
\begin{enumerate}[leftmargin=2em]
\item
whether setup and solve learning reduce end-to-end runtime;
\item how performance evolves
during online learning;
\item how the method
behaves on a more difficult variable-flow problem family;
\item whether adaptive relaxation provides
benefit when the hierarchy is held fixed.
\end{enumerate}
We first compare BoomerAMG's default, LSTDQ, LinUCB, and LinUCB--LSTDQ across problem sizes and random seeds, then examine the online trajectories and cycle-level relaxation choices. A matched-hierarchy study isolates the effect of relaxation control from setup selection, and a diffusion--advection test with a spatially varying flow field evaluates robustness beyond the constant-coefficient sequences. The primary metric is total recorded runtime, including recurring learner computation; setup-plus-solve execution time is also reported separately. The code and experiment scripts are publicly available\footnote{\url{https://github.com/EriseHe/BoomerAMG-Autotuner}}.

All runs use the same source and solver-library versions on an Apple M5 Pro with 24 GB RAM, running macOS 27.0 (ARM64), with HYPRE 3.0.0 and Python 3.10.21, one MPI rank, and one compute thread per solver process. The four methods execute serially within each problem using balanced random-permutation ordering.

\subsection{Problem families, features, and protocol}
\label{subsec:experiment_setup}
We compare four complete online procedures. 
Throughout the experiments, Default denotes BoomerAMG's default setup configuration with fixed relaxation weight \(w=1\) and the common \(\ell_1\)-Jacobi smoothing profile specified below; it uses one attempt per problem, \(J_{\max}=1\). LSTDQ keeps the default setup but learns the cycle-level relaxation weights; if its learned solve attempt fails, one default fallback is allowed with \(J_{\max}=2\). LinUCB learns the setup parameters while retaining the reference solve settings, and the combined LinUCB--LSTDQ learns both the setup and solve. The two setup-learning methods, namely LinUCB and LinUCB--LSTDQ, allow up to three learned setup attempts followed by one default fallback where \(J_{\max}=4\). The reported runtime includes all executed attempts.
Each learned method maintains and updates its own model. Thus, the four-way comparison evaluates complete online procedures rather than different relaxation policies on the same hierarchy. Section~\ref{subsec:matched_policy} isolates the relaxation effect using matched hierarchies.

\paragraph{Diffusion family}
The first problem family is the seven-point finite-difference discretization of
\begin{equation}
    -\partial_x(c_x\partial_x u)
    -\partial_y(c_y\partial_y u)
    -\partial_z(c_z\partial_z u)=f
    \quad\text{in }(0,1)^3,
    \qquad
    u|_{\partial(0,1)^3}=0.
    \label{eq:diffusion_pde}
\end{equation}
The diffusion--convection problem generator in \hypre (DifConv) is used on 
the grids of sizes \(40^3\), \(60^3\), and \(80^3\), 
giving \(64{,}000\), \(216{,}000\), and \(512{,}000\) unknowns. For each problem, \(c_x,c_y,c_z\) are sampled independently and uniformly from \([1,1000]\) and remain spatially constant; advection is zero.
The learners use normalized PDE coefficients as context, rather than spectral information extracted from the matrix. Define
\begin{equation}
    \ell_{i,t}
    =
    \frac{\log(c_{i,t})}{\log(1000)},
    \qquad i\in\{x,y,z\}.
\end{equation}
Both learners receive
\begin{equation}
    \chi_t
    =
    (1,\ell_{x,t},\ell_{y,t},\ell_{z,t})^\top
    \in\mathbb R^4.
    \label{eq:experiment_context}
\end{equation}

\paragraph{Diffusion--advection family}
The second family discretizes
\begin{equation}
    -\nabla\!\cdot(\mathbf D_t\nabla u)
    +\mathbf b_t\cdot\nabla u=f,
    \qquad
    \mathbf D_t
    =
    \operatorname{diag}(c_{x,t},c_{y,t},c_{z,t}),
    \label{eq:advection_pde}
\end{equation}
on the same domain and grids, with homogeneous Dirichlet boundary conditions. The generator uses centered differences for diffusion and upwind differences for advection. The three diffusion and three advection coefficients are sampled independently and uniformly from \([1,1000]\) and remain spatially constant within each problem.
Using the same \(\ell_{i,t}\) as above and \(b_{\rm scale}=1000\), define
\begin{equation}
    \widetilde b_{i,t}
    =
    \frac{
        \operatorname{sgn}(b_{i,t})
        \log(1+|b_{i,t}|)
    }{
        \log(1+b_{\rm scale})
    },
    \qquad i\in\{x,y,z\}.
\end{equation}
The context is
\begin{equation}
    \chi_t
    =
    (1,\ell_{x,t},\ell_{y,t},\ell_{z,t},
       \widetilde b_{x,t},\widetilde b_{y,t},
       \widetilde b_{z,t})^\top
    \in\mathbb R^7.
    \label{eq:advection_context}
\end{equation}

\paragraph{Solver and online protocol}
All experiments use seeded random right-hand sides with zero initial iterates. The common solve profile uses one \(\ell_1\)-Jacobi pre- and post-smoothing sweep and Gaussian elimination on the coarsest level, which correspond to \hypre relaxation options 18 and 9. The reference relaxation weight is \(w=1\).
Every attempt uses relative residual tolerance \(10^{-6}\)
and cycle limit \(K_{\max}=50\). %

At each grid size, all four methods receive the same ordered sequence of 5000 matrix/right-hand-side pairs. The learned methods start without pretraining and update throughout this online phase. They then process 1000 additional problems with model updates and random exploration disabled. Recurring overhead, including feature construction and action selection, remains included in both phases. Both of the setup learners use the search space in Table~\ref{tab:search_space} and the same 1024-candidate construction procedure, but evolve independently. Each solve learner likewise receives the setup descriptor \(\zeta(a)\) selected by its own branch.
We report \ModuleFourSeedCoverage{}. Each seed determines the problem sequence, learner randomization, candidate generation, and within-problem method ordering.

\paragraph{Timing and aggregation}
The solver cost \(C_t^{\rm solver}\) in \eqref{eq:timing_decomposition} includes setup and solve execution across all attempts. The total cost \(C_t^{\rm total}\) additionally includes recurring setup- and solve-learner computation. The times of matrix and right-hand-side assembly, experiment orchestration, one-time preparation, and output I/O are excluded. Failed attempts remain included in both costs.

We report runtime reductions over the full online learning sequence and the additional frozen phase to evaluate generalization to new instances. The corresponding index sets are
\begin{equation}
    W_{\rm online}=\{1,\ldots,5000\},
    \qquad W_{\rm frozen}=\{5001,\ldots,6000\}.
\end{equation}
For seed \(s\), window \(W\in \{W_{\rm online},W_{\rm frozen}\}\), and method \(\Pi\), the total-time reduction relative to the default is
\begin{equation}
    R^{\rm total}_{s,W}(\Pi)
    =
    100\left[
        1-
        \frac{
            \sum_{t\in W}C_{s,t}^{\rm total}(\Pi)
        }{
            \sum_{t\in W}C_{s,t}^{\rm total}(\Pi_{\rm def})
        }
    \right].
    \label{eq:reported_reduction}
\end{equation}
All problems ultimately succeed in the constant-coefficient experiments, so these sums include the full cost of any failed attempts and subsequent recovery. Solver-time reduction is defined analogously using \(C^{\rm solver}\).
Reported reductions are the mean and sample standard deviation of paired seed-level values, with seeds weighted equally. The standard deviation therefore measures variation across independent adaptive runs; individual problems within a sequence are not treated as independent replicates.

\subsection{Runtime across problem sizes}
\label{subsec:scale_results}

Table~\ref{tab:scale_runtime} summarizes online runtime and total-time reductions for all four methods at the three grid sizes. Setup and solve execution times are reported separately from learner overhead, and total-time reductions include all three components. Across the \(40^3\), \(60^3\), and \(80^3\) grids, LinUCB--LSTDQ reduces total runtime by \ModuleFourDiffusionFullSavings{} for diffusion and \ModuleFourAdvectionFullSavings{} for diffusion--advection. Its learner overhead accounts for approximately 11.8--13.2\%, 4.2--4.7\%, and 1.8--2.0\% of total online runtime, respectively.

\begin{table}[t]
\centering
\caption{Online runtime decomposition and total-time reductions during the online (problems 1--5000) and frozen (problems 5001--6000) phases. Setup and solve execution include all attempts, and learner overhead includes recurring setup- and solve-learner computation. Time columns are in ms per problem and average the six seed-level means. Both reduction columns report the mean \(\pm\) sample standard deviation of paired total-time reductions relative to the default, including learner overhead; positive values favor the learned method.}
\label{tab:scale_runtime}
\footnotesize
\setlength{\tabcolsep}{4pt}
\newlength{\runtimecostwidth}
\newlength{\runtimeoverheadwidth}
\newlength{\runtimereductionwidth}
\settowidth{\runtimecostwidth}{000.0}
\settowidth{\runtimeoverheadwidth}{00.0}
\settowidth{\runtimereductionwidth}{\(00.0\pm0.0\)\%}
\newcolumntype{R}[1]{>{\hbox to #1\bgroup\hss}c<{\egroup}}
\begin{tabular}{@{}lR{\runtimecostwidth}R{\runtimecostwidth}R{\runtimeoverheadwidth}R{\runtimereductionwidth}R{\runtimereductionwidth}@{}}
\toprule
& \multicolumn{4}{c}{Online (1--5000)} & \multicolumn{1}{c}{Frozen (5001--6000)} \\
\cmidrule(lr){2-5}\cmidrule(l){6-6}
Method & \multicolumn{1}{c}{Setup} & \multicolumn{1}{c}{Solve}
& \multicolumn{1}{c}{\shortstack{Learner\\overhead}}
& \multicolumn{1}{c}{\shortstack{Total-time\\reduction}}
& \multicolumn{1}{c}{\shortstack{Total-time\\reduction}} \\
\midrule
\multicolumn{6}{@{}l}{\textit{Diffusion, \(40^3\)}} \\
The default & 57.4 & 24.9 & 0.0 & --- & --- \\
LSTDQ & 57.4 & 17.5 & 2.0 & \(6.5\pm0.5\)\% & \(8.1\pm0.7\)\% \\
LinUCB & 21.9 & 40.7 & 3.1 & \(20.2\pm0.7\)\% & \(28.4\pm1.3\)\% \\
LinUCB--LSTDQ & 20.7 & 25.6 & 6.2 & \(36.1\pm1.0\)\% & \(47.0\pm1.7\)\% \\
\addlinespace[3pt]
\multicolumn{6}{@{}l}{\textit{Diffusion, \(60^3\)}} \\
The default & 203.7 & 88.4 & 0.0 & --- & --- \\
LSTDQ & 204.2 & 64.6 & 2.5 & \(7.1\pm0.6\)\% & \(9.2\pm0.6\)\% \\
LinUCB & 73.6 & 135.7 & 3.1 & \(27.3\pm1.1\)\% & \(32.9\pm0.9\)\% \\
LinUCB--LSTDQ & 70.7 & 87.9 & 6.9 & \(43.3\pm1.3\)\% & \(49.8\pm0.8\)\% \\
\addlinespace[3pt]
\multicolumn{6}{@{}l}{\textit{Diffusion, \(80^3\)}} \\
The default & 514.3 & 232.4 & 0.0 & --- & --- \\
LSTDQ & 515.3 & 161.4 & 2.8 & \(9.0\pm0.4\)\% & \(9.8\pm0.4\)\% \\
LinUCB & 208.8 & 350.8 & 3.8 & \(24.6\pm4.0\)\% & \(29.9\pm2.4\)\% \\
LinUCB--LSTDQ & 178.9 & 242.3 & 7.9 & \(42.5\pm2.3\)\% & \(49.1\pm1.3\)\% \\
\addlinespace[3pt]
\multicolumn{6}{@{}l}{\textit{Diffusion--advection, \(40^3\)}} \\
The default & 58.3 & 24.6 & 0.0 & --- & --- \\
LSTDQ & 58.3 & 17.3 & 2.2 & \(6.2\pm0.4\)\% & \(7.9\pm0.5\)\% \\
LinUCB & 28.7 & 42.6 & 4.2 & \(8.9\pm2.7\)\% & \(24.8\pm3.2\)\% \\
LinUCB--LSTDQ & 23.0 & 34.7 & 8.8 & \(19.8\pm1.3\)\% & \(37.0\pm1.7\)\% \\
\addlinespace[3pt]
\multicolumn{6}{@{}l}{\textit{Diffusion--advection, \(60^3\)}} \\
The default & 206.9 & 87.9 & 0.0 & --- & --- \\
LSTDQ & 207.4 & 62.2 & 2.7 & \(7.6\pm0.7\)\% & \(9.1\pm0.7\)\% \\
LinUCB & 100.0 & 147.9 & 4.3 & \(14.5\pm2.5\)\% & \(25.6\pm2.2\)\% \\
LinUCB--LSTDQ & 78.1 & 115.8 & 9.6 & \(31.0\pm2.0\)\% & \(40.5\pm2.1\)\% \\
\addlinespace[3pt]
\multicolumn{6}{@{}l}{\textit{Diffusion--advection, \(80^3\)}} \\
The default & 521.1 & 229.6 & 0.0 & --- & --- \\
LSTDQ & 522.0 & 159.2 & 3.0 & \(8.9\pm0.3\)\% & \(9.8\pm0.3\)\% \\
LinUCB & 222.8 & 413.9 & 4.5 & \(14.6\pm1.6\)\% & \(23.1\pm2.7\)\% \\
LinUCB--LSTDQ & 201.8 & 311.4 & 10.6 & \(30.2\pm1.6\)\% & \(36.8\pm1.4\)\% \\

\bottomrule
\end{tabular}
\end{table}

Table~\ref{tab:scale_runtime} also compares the full online sequence and frozen phase. For LinUCB--LSTDQ, the reductions in the additional frozen phase range from \ModuleFourDiffusionFrozenSavings{} for diffusion and from \ModuleFourAdvectionFrozenSavings{} for diffusion--advection. The gains are smaller for diffusion--advection than for diffusion, but the frozen-phase reductions exceed the full online-sequence reductions at every grid size. This pattern suggests that the policy learned over the first 5,000 problems remains effective on new instances.

\subsection{Runtime decomposition and online behavior}
\label{subsec:runtime_results}
\begin{figure}[!htbp]
\centering
\includegraphics[width=0.8\linewidth]{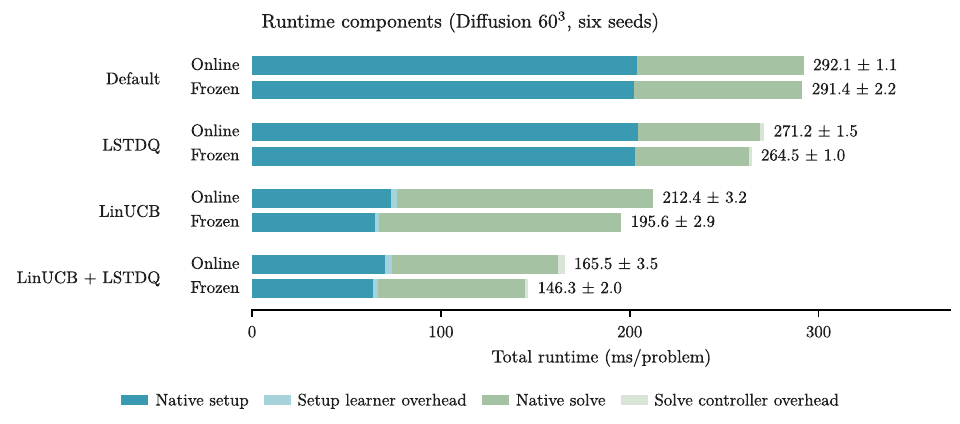}
\par\nointerlineskip
\setlength{\abovecaptionskip}{6pt}
\caption{Runtime decomposition on \(\ModuleFourGrid^3\) diffusion across \(\ModuleFourSeedCountWord{}\) seeds. Paired stacks show setup execution, setup-learner overhead, solve execution, and solve-learner overhead. Labels report total runtime as mean \(\pm\) sample standard deviation of seed-level means (ms/problem), for the online phase (problems 1--5000) and the frozen phase (problems 5001--6000).}
\label{fig:runtime_breakdown}
\end{figure}

Figure~\ref{fig:runtime_breakdown} shows how setup, solve, and learner costs contribute to total runtime on the \(\ModuleFourGrid^3\) diffusion problems. The online bars summarize the full sequence, problems \(1\)--\(5000\). Setup-only LinUCB reduces mean setup time by 63.9\% per problem, but increases mean solve time by 53.5\% relative to the default. Thus, a hierarchy that is cheaper to construct may not yield a more efficient solve phase.
Standalone LSTDQ reduces mean solve time by 26.9\% on the default hierarchy (88.4 to 64.6 ms). Across the independently trained branches, LinUCB--LSTDQ has 35.2\% lower mean solve time than setup-only LinUCB (87.9 versus 135.7 ms), with similar setup costs (70.7 versus 73.6 ms per problem). The joint method therefore preserves the setup savings while reducing the solve-cost penalty observed with setup-only learning.
The frozen bars, which report problems \(5001\)--\(6000\), show the same qualitative pattern on new instances after the learners have accumulated experience. Together, the paired bars show that setup selection and cycle-level relaxation control address complementary components of the total runtime.

\begin{figure}[!htbp]
\centering
\includegraphics[width=\linewidth]{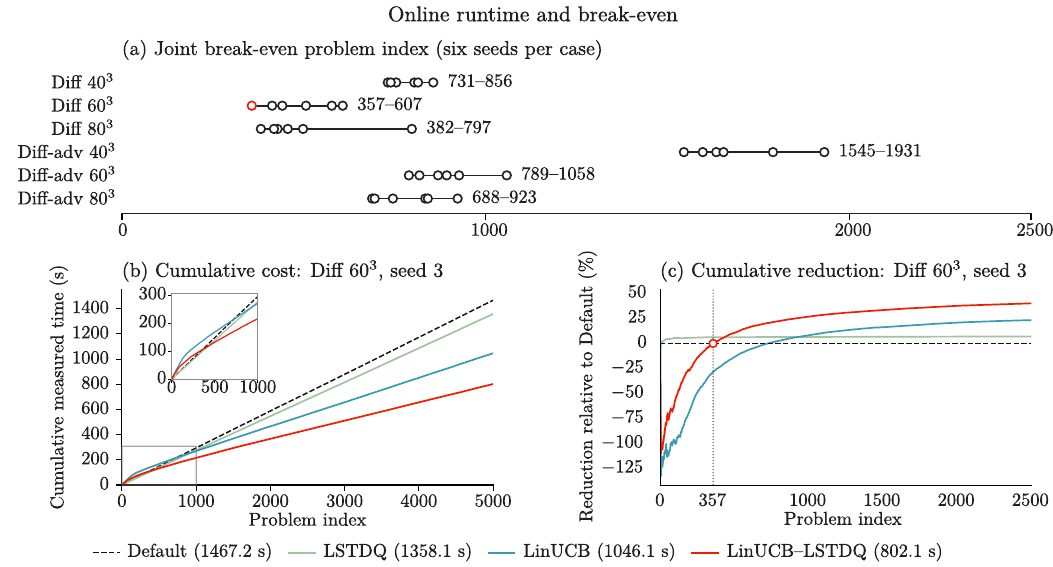}
\par\nointerlineskip
\setlength{\abovecaptionskip}{6pt}
\caption{Online runtime and break-even behavior. (a) Joint-method break-even problem indices across six seeds per case. The break-even problem index is the first index from which the joint method's cumulative total runtime remains below the default's through problem 5000. Open points show individual seeds; whiskers and labels give the minimum--maximum range. The red point identifies the run illustrated in (b,c). (b,c) Illustrative \(\ModuleFourGrid^3\) diffusion run, seed \ModuleFourSeed{}: cumulative total runtime and percentage reduction relative to the default, over problems 1--5000 in (b) and 1--2500 in (c). All executed attempts and recurring learner computation are included. Curves are unsmoothed prefix sums; the inset enlarges problems 1--1000. The black vertical guide in (c) marks the break-even problem index, 357. Legend values give online-phase runtime in seconds.}
\label{fig:cumulative}
\end{figure}

Across all three grid sizes and learned methods, the cumulative runtime reduction increases over the 5,000-problem sequence. To illustrate how the performance difference develops during online learning across problem families and grid sizes, we show the first problem index from which LinUCB–LSTDQ overcomes its initial runtime disadvantage and maintains cumulative savings over the default for the remainder of the sequence as shown in Figure~\ref{fig:cumulative}. Smaller PDE systems, such as those on a $40^3$ grid, have shorter runtimes, making learner overhead a larger fraction of the total cost and thus harder to amortize. Diffusion–advection problems generally reach payback later than pure diffusion problems. In particular, we show one \(\ModuleFourGrid^3\) diffusion trajectory to illustrate how the runtime advantage develops progressively over the sequence rather than arising from a small number of isolated problems. The selected seed, \ModuleFourSeed{}, has the largest full-sequence reduction among the completed diffusion runs and is used only for illustration; the aggregate comparisons above provide the primary performance results.
For this run, LinUCB--LSTDQ reduces average total runtime from \ModuleFourSelectedDefaultCost{} to \ModuleFourSelectedJointCost{} ms per problem over the full sequence, a \ModuleFourSelectedFullReduction{}\% reduction relative to the default. Over the 1000 frozen problems, the reduction is \ModuleFourSelectedFrozenReduction{}\%.

\subsection{Hierarchy and cycle-level parameter choices}
\label{subsec:sequential_results}

\begin{figure}[!htbp]
\centering
\includegraphics[width=\linewidth]{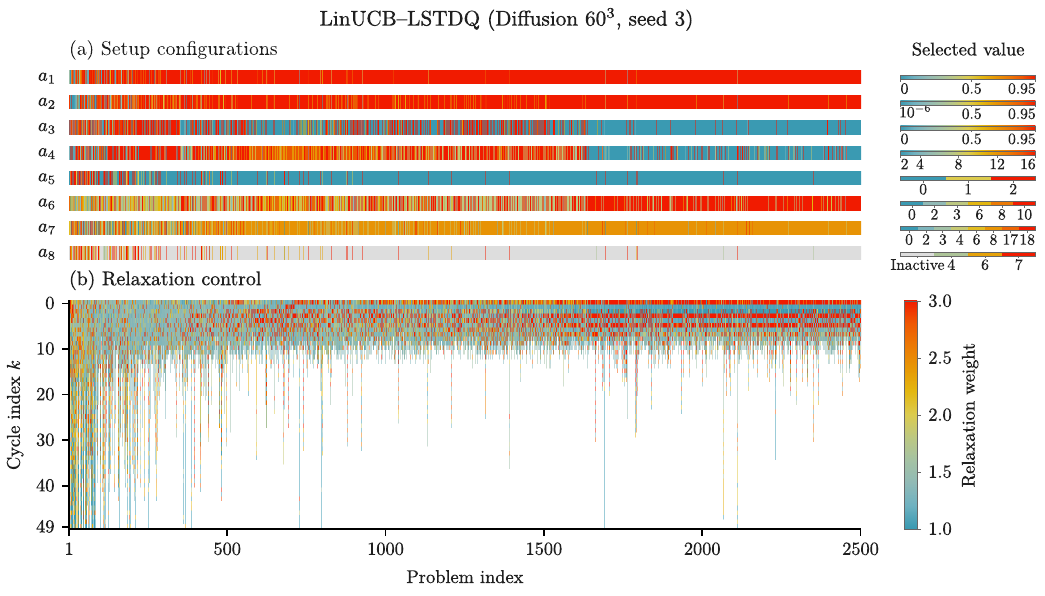}
\par\nointerlineskip
\setlength{\abovecaptionskip}{6pt}
\caption{First-attempt setup choices and relaxation weights selected by LinUCB--LSTDQ in the illustrative \(\ModuleFourGrid^3\) diffusion run, seed \ModuleFourSeed{}. Columns correspond to problems 1--2500. (a) Setup parameters \(a_1,\ldots,a_8\), defined in Table~\ref{tab:search_space}. (b) Rows correspond to zero-based cycle indices, increasing downward. White cells indicate that the first attempt has terminated. Cycles from retries and the default fallback are omitted.}
\label{fig:weights}
\end{figure}

In the selected 5000-problem online run, LinUCB--LSTDQ has 88 first-attempt failures, all of which are resolved by retries without using the default fallback. All failed attempts and subsequent recovery are included in the reported runtime. Figure~\ref{fig:weights} shows how first-attempt setup configurations and relaxation weights evolve as the learners train online from scratch in this run.
In Figure~\ref{fig:weights}(a), the setup choices become increasingly concentrated as more problems are solved. Over problems 2001--2500, first attempts favor a strength threshold of 0.95 (in 99.6\% of the problems), standard interpolation (97.4\%), and a maximum of two nonzero entries per interpolation row (83.0\%). The higher strength threshold restricts the strong-connection graph, while the smaller row limit restricts interpolation density. These choices can reduce construction cost, consistent with the benefits of limiting interpolation complexity reported for three-dimensional diffusion problems~\cite{desterck2008distance}. The learned relaxation policy also develops more prominent recurring high--low patterns as the hierarchy selections narrow. As shown in Figure~\ref{fig:weights}(b), its choices vary across both problems and cycle indices. Table~\ref{tab:best_individual_configuration} gives a concrete example from problem 2275 in this run. Benefiting from both more efficient hierarchy construction and a nonstationary relaxation policy, the method achieves the largest recorded individual total-time reduction, 72.7\%, relative to the default. Despite higher grid complexity, the operator and memory complexities are 15.4\% and 18.2\% lower, respectively, which helps explain the observed 81.0\% reduction in setup time. The 51.7\% reduction in native solve time results from both fewer V-cycles (13 to 8) and a 21.5\% reduction in average time per V-cycle (7.71 to 6.05 ms).

\begin{table}[!htbp]
\centering
\caption{Setup configuration, hierarchy statistics, and runtime for the \(60^3\) diffusion problem 2275, seed 3. Only setup parameters that differ from the default are shown. Both methods converge on their first attempt. Total time includes learner overhead; the per-problem total-time reduction is 72.7\%. This selected best case is illustrative rather than an average performance result. Grid complexity is the ratio of total unknowns across levels to fine-grid unknowns. Operator and memory complexities are the ratios of total operator nonzeros and total operator-plus-interpolation nonzeros, respectively, to fine-grid operator nonzeros.}
\label{tab:best_individual_configuration}
\begingroup
\footnotesize
\setlength{\tabcolsep}{4pt}
\begin{tabularx}{\linewidth}{@{}cXccc@{}}
\toprule
Symbol & Parameter & Default & LinUCB--LSTDQ & Reduction \\
\midrule
\(a_1\) & \texttt{strong\_threshold} & 0.25 & 0.95 & \\
\(a_2\) & \texttt{max\_row\_sum} & 0.90 & 0.95 & \\
\(a_4\) & \texttt{P\_max\_elmts} & 4 & 2 & \\
\(a_7\) & \texttt{interp\_type} & Extended+i (6) & Standard (8) & \\
\midrule
\multicolumn{2}{@{}l}{Setup execution (ms)} & 308.84 & 58.54 & 81.0\% \\
\multicolumn{2}{@{}l}{\hspace{1em}\textit{Grid complexity}} & 1.676759 & 1.971565 & $-17.6$\% \\
\multicolumn{2}{@{}l}{\hspace{1em}\textit{Operator complexity}} & 4.493422 & 3.800977 & 15.4\% \\
\multicolumn{2}{@{}l}{\hspace{1em}\textit{Memory complexity}} & 5.166080 & 4.224497 & 18.2\% \\
\addlinespace[2pt]
\multicolumn{2}{@{}l}{Solve execution (ms)} & 100.20 & 48.41 & 51.7\% \\
\multicolumn{2}{@{}l}{\hspace{1em}\textit{V-cycles}} & 13 & 8 & 38.5\% \\
\multicolumn{2}{@{}l}{\hspace{1em}\textit{Average time per V-cycle (ms)}} & 7.71 & 6.05 & 21.5\% \\
\addlinespace[2pt]
\multicolumn{2}{@{}l}{Setup- and solve-learner overhead (ms)} & 0 & 4.58 & --- \\
\multicolumn{2}{@{}l}{Total time (ms)} & 409.03 & 111.53 & 72.7\% \\
\bottomrule
\end{tabularx}
\endgroup
\end{table}

The combination of the learned AMG hierarchy and solve policy is consistent with the complementary roles of hierarchy selection and relaxation control. Coarse-grid correction and relaxation divide the work of reducing error. Accordingly, changing interpolation and coarse-grid operators during setup changes which error components coarse-grid correction handles effectively---and hence which components relaxation must reduce~\cite{xu2017algebraic,huang2023learning}. Section~\ref{subsec:matched_policy} provides a numerical rationale for the recurring high--low relaxation weights observed in Figure~\ref{fig:weights}(b) and compares fixed, prescribed periodic, and learned policies on identical hierarchies. Over the full displayed sequence, problems 1--2500, mean setup time decreases from 203.1 to 76.1 ms per problem, while mean solve time remains close to that of the default (92.3 versus 90.0 ms). Including all attempts and learner overhead, total runtime decreases by 40.1\%.

\begin{figure}[!htbp]
\centering
\includegraphics[width=\linewidth]{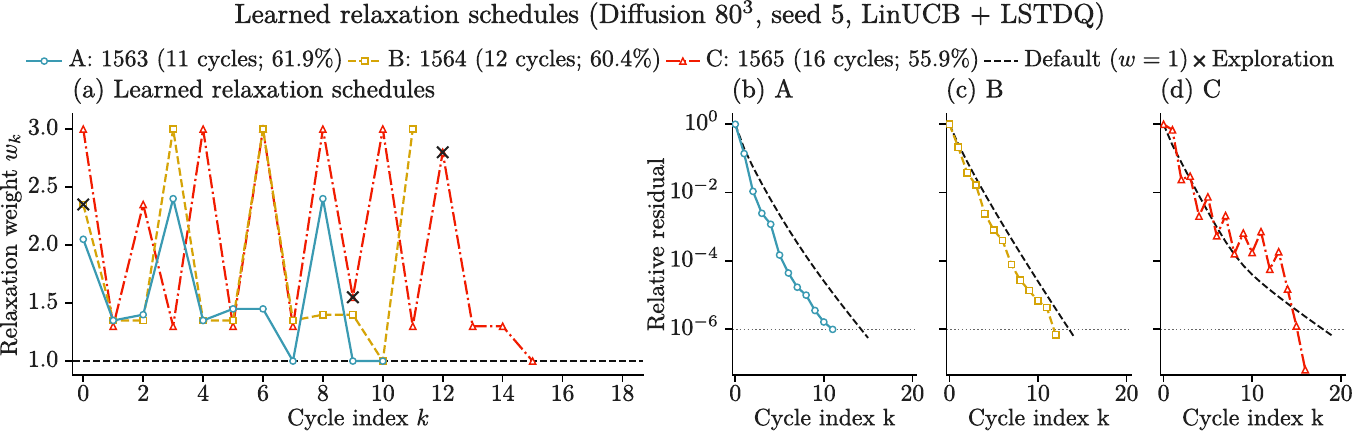}
\par\nointerlineskip
\setlength{\abovecaptionskip}{6pt}
\caption{Adaptive relaxation on \(\ModuleFourExampleGrid^3\) diffusion, seed \ModuleFourExampleSeed{}, for consecutive problems A, B, and C (\ModuleFourExampleIDs{}). (a) Recorded LinUCB--LSTDQ weights; crosses mark random exploratory actions, and the horizontal black line marks \(w=1\). (b--d) Residual histories for problems A--C, respectively, comparing LinUCB--LSTDQ (colored) with the paired default (black). The three residual panels share axis scales. Dotted lines mark the convergence tolerance. Legend percentages give total-time reductions relative to the default on the same problem, including learner overhead. All three joint-method solves succeed on their first attempt and are shown through convergence.}
\label{fig:individual_cycle_schedules}
\end{figure}

To examine this behavior at the level of individual solves, Figure~\ref{fig:individual_cycle_schedules} shows three consecutive problems A, B, and C (\ModuleFourExampleIDs{}) from the joint method on \(\ModuleFourExampleGrid^3\) diffusion, seed \ModuleFourExampleSeed{}. Under the online-selected setup configurations, the three problems reach tolerance in \ModuleFourExampleCycles{} cycles with different relaxation trajectories. Their initial relaxation weights are \(\ModuleFourExampleInitialWeights{}\), respectively, and only the initial choice for problem B is a random exploratory action. The corresponding residual histories in panels (b)--(d) show how the solves progress. For reference, Figure~\ref{fig:individual_cycle_schedules} also shows the default on the same problem inputs, using the default hierarchy and fixed weight \(w=1\). In particular, problem C exhibits temporary residual increases during several high-weight cycles, followed by substantial decreases under lower weights. Despite this nonmonotone behavior, the solve reaches tolerance in 16 cycles, compared with 19 for the default; the joint method reduces total runtime by 55.9\% on that problem. These examples illustrate the cycle-level behavior produced by the online learner, but they do not by themselves isolate the effect of relaxation control because the joint method and the default use different hierarchies. The matched-hierarchy experiment in Section~\ref{subsec:matched_policy} provides that controlled comparison.

\subsection{Variable-flow diffusion--advection}
\label{subsec:variable_flow}
Next, we evaluate the online methods on a more challenging diffusion--advection family with a spatially varying flow field. This experiment complements the preceding constant-coefficient tests by introducing stronger variation in the linear systems encountered along the online sequence. We use \hypre's generator\footnote{\url{https://github.com/hypre-space/hypre/blob/v3.0.0/src/parcsr_ls/par_vardifconv_rs.c}} to discretize the following problem:
\begin{equation}
    \begin{gathered}
    -\epsilon_t\Delta u
    +\mathbf v(x,y,z)\cdot\nabla u=f
    \quad\text{in }(0,1)^3,
    \qquad
    u|_{\partial(0,1)^3}=0,\\
    \mathbf v(x,y,z)
    =
    \bigl(
        (2y-1)(1-x^2),
        2xy(y-1),
        2xy(y-1)
    \bigr).
    \end{gathered}
    \label{eq:variable_flow_pde}
\end{equation}
The flow field is fixed, while \(\log_{10}\epsilon_t\) is sampled uniformly for each problem from \([-4,-1]\), so the diffusion strength varies over three orders of magnitude. We use a \(60^3\) interior grid, seeded random right-hand sides, and zero initial iterates. The learner context is
\begin{equation}
    \chi_t
    =
    \left(
        1,\,
        \frac{2(\log_{10}\epsilon_t+4)}{3}-1
    \right)^\top .
\end{equation}
As in the constant-coefficient experiments, the learned methods begin without pretraining and learn online from the first problem in each sequence. We use six paired 5000-problem sequences. The setup search space, smoother, stopping tolerance, cycle limit, and method-specific attempt budgets are unchanged from Section~\ref{subsec:experiment_setup}. Because the methods have different recovery budgets, we report both first-attempt failures and failures that remain after each permitted attempt.

\begin{table}[!htbp]
\centering
\caption{Variable-flow diffusion--advection on the \(60^3\) grid, over 30,000 problems across six seeds. Failure columns give the numbers of problems remaining unsolved after each attempt; italicized counts mark final failures. The last permitted attempt uses the default recovery for LSTDQ, LinUCB, and LinUCB--LSTDQ. The final column gives total-time reductions relative to the paired default as mean \(\pm\) sample standard deviation across six seeds, including all executed attempts and recurring learner overhead.}
\label{tab:variable_flow}
\footnotesize
\setlength{\tabcolsep}{4pt}
\begin{tabular}{@{}lllllr@{}}
\toprule
\multirow{2}{*}[-0.3ex]{Method} & \multicolumn{4}{c}{Failures after attempt}
& \multirow{2}{*}[-0.3ex]{\shortstack{Total-time\\reduction}} \\
\cmidrule(lr){2-5}
& {\footnotesize \(j=1\)} & {\footnotesize \(j=2\)}
& {\footnotesize \(j=3\)} & {\footnotesize \(j=4\)} & \\
\midrule
Default (\(J_{\max}=1\)) & \textit{1677} & --- & --- & --- & --- \\
LSTDQ (\(J_{\max}=2\)) & 777 & \textit{26} & --- & --- & \(9.8\pm5.6\%\) \\
LinUCB (\(J_{\max}=4\)) & 810 & 62 & 44 & \textit{28} & \(50.8\pm2.9\%\) \\
LinUCB--LSTDQ (\(J_{\max}=4\)) & 1152 & 138 & 69 & \textit{2} & \(54.1\pm2.1\%\) \\

\bottomrule
\end{tabular}
\end{table}

The differences between the methods are substantially larger on this variable-flow family. LinUCB--LSTDQ reduces the number of final failures from \(1677\) for the default to only \(2\), while reducing total time by \(54.1\pm2.1\%\) over the full sequence. Over the final 1000 problems of each sequence, the joint method has one final failure among 6000 problems and achieves a \(59.2\pm2.4\%\) total-time reduction. The stronger late-sequence result is consistent with the online learners becoming more effective as additional problems are observed.

The runtime reduction is not explained simply by the methods succeeding on different sets of problems. Restricting the comparison to the \(28,323\) problems completed by all four methods gives a similar full-sequence total-time reduction of \(53.4\pm2.2\%\) for LinUCB--LSTDQ.
The fixed-setup comparison isolates the contribution of cycle-level relaxation control. With the reference hierarchy held fixed, the learned LSTDQ attempt succeeds without fallback on \(1,651\) of the \(1,677\) problems on which Default fails, within the same 50-cycle limit. Including recovery attempts, its mean solve time decreases from \(304.3\) ms per problem for the default to \(230.3\) ms. Thus, adaptive relaxation alone provides both a substantial robustness improvement and a reduction in solve cost on this problem family.
The joint method obtains additional gains by adapting the setup configuration. Including recovery attempts, its mean setup time is \(129.8\) ms per problem, compared with \(324.6\) ms for the default, while its mean solve time decreases further to \(150.3\) ms. These phase-level timings show that the overall improvement comes from both components: \emph{Standalone LSTDQ improves solve performance on a fixed hierarchy, while the joint method uses LinUCB to reduce setup cost and select hierarchies that further benefit adaptive relaxation.}

\subsection{Matched-hierarchy comparison of relaxation policies}
\label{subsec:matched_policy}

We then perform matched-hierarchy experiments in which the hierarchy is held fixed and only the relaxation policy is changed. The experiments ask two questions:
\begin{enumerate}[leftmargin=2em]
\item whether varying
relaxation weights across cycles improves over the fixed-weight strategy;
\item whether the learned LSTDQ policy improves
over a prescribed varying-weight schedule derived from classical smoothing analysis.
\end{enumerate}
We use the final learner states from the six \(60^3\) diffusion and six \(60^3\) variable-flow diffusion--advection online runs and freeze them for this controlled evaluation. Each frozen learner is evaluated on the same 100 additional problems that were not part of its original 5000-problem sequence. For each learner state and evaluation problem, the setup configuration selected by the frozen setup learner is held fixed across all relaxation policies. Thus, the compared policies use the same matrix, right-hand side, zero initial iterate, hierarchy, tolerance, and cycle limit. No evaluation observation updates either learner.
Results are then aggregated over the 100 problems, and variability is reported across the six independent online runs.
For each frozen learner state, we compare five policies:
\begin{enumerate}[leftmargin=2em]
    \item \emph{Fixed \(w\):} the reference relaxation
    weight used by the default.

    \item \emph{Global best:} a single best fixed weight
    chosen from the 41 admissible weights to minimize total first-attempt native solve time over the 100 evaluation problems.

    \item \emph{Per-problem best:} the best fixed
weight chosen separately from the 41 admissible weights to minimize first-attempt native solve time for each evaluation problem.

    \item \emph{Period two:} a prescribed alternating
two-weight schedule derived from an SPD smoothing surrogate, without using measured PDE runtimes.

    \item \emph{Frozen LSTDQ:} the cycle-level policy from
the corresponding completed online run, with its model parameters fixed during evaluation.
\end{enumerate}
The global and per-problem best fixed policies are comparison references rather than deployable online methods, since their selection requires evaluating all 41 candidate weights on the evaluation problems. The weight-selection scan averages three timings per problem--weight pair. After selection, the chosen policies are retimed independently; the weight-selection scan itself is excluded from the reported cost.

The period-two policy provides a prescribed nonlearning reference based on the weighted two-level smoothing criterion of Lottes~\cite{lottes2023optimal}. Its construction exploits the SPD structure of the diffusion problems rather than measured runtimes. For an SPD matrix \(\mathbf A=(a_{ij})\), let \(D=\operatorname{diag}(\sum_j|a_{ij}|)\). Then \(\mathbf A\preceq D\). The eigenvalues of \(D^{-1}\mathbf A\) lie in \((0,1]\). For weighted Jacobi, the error update is \(e'=(I-wD^{-1}\mathbf A)e\). For an eigenmode with eigenvalue \(\xi\), one step multiplies its error by \(1-w\xi\). For the smoothing surrogate, we consider two successive weighted-Jacobi steps with weights \(a\) and \(b\). These steps multiply the mode's error by
\begin{equation}
    p_{a,b}(\xi)=(1-a\xi)(1-b\xi),
\end{equation}
where a small \(|p_{a,b}(\xi)|\) means that the mode is strongly reduced. The quantity
\begin{equation}
    \eta(a,b)
    =
    \max_{0\le\xi\le1}
    \xi\,p_{a,b}(\xi)^2
    \label{eq:period_two_objective}
\end{equation}
weights the squared error amplification by \(\xi\) and takes the worst case over the spectral interval. This emphasizes modes with larger eigenvalues, which relaxation is intended to damp. Modes near \(\xi=0\) receive less weight because coarse-grid correction is expected to handle them. Minimizing \(\eta(a,b)\) gives a prescribed alternating pair for comparison with learned weights. Note here that this is a smoothing surrogate, rather than a prediction of full AMG convergence or runtime. The unconstrained optimum is \(2\pm2/\sqrt5\). Because the solve learner only permits weights in the discrete action set \(\mathcal W\), we evaluate the surrogate over all pairs \(a,b\in\mathcal W\):
\begin{equation}
    (a_{\mathcal W},b_{\mathcal W})=(2.85,1.10)
    \in
    \operatorname*{arg\,min}_{a,b\in\mathcal W}
    \eta(a,b).
    \label{eq:period_two_grid}
\end{equation}
We apply the resulting pair in high--low order between complete AMG cycles,
\begin{equation}
    2.85,\ 1.10,\ 2.85,\ 1.10,\ldots,
\end{equation}
with the sequence restarted at the beginning of each problem.

This construction should be interpreted as a strong structured baseline for the present SPD diffusion test. In particular, the same derivation does not directly apply to the nonsymmetric diffusion--advection problems considered elsewhere in this paper.
All five policies are timed under the same evaluation protocol. We compare the native solver execution time in the first attempt on the matched hierarchy. Hierarchy-construction time, learner computation, and all subsequent recovery or fallback work are excluded. Elapsed times from all 100 problems are retained regardless of convergence; these timings do not measure time to successful completion. Failure rates are reported alongside them and refer to the prescribed \(K_{\text{max}}=50\) stopping rule.

\begin{figure}[!htbp]
\centering
\includegraphics[width=\linewidth]{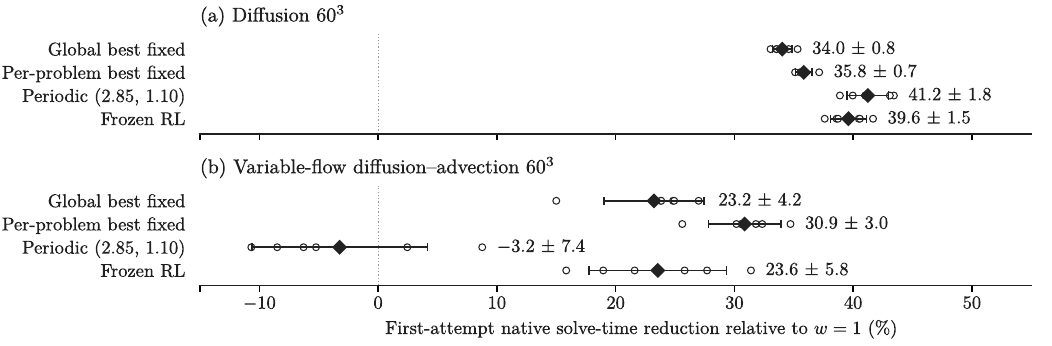}
\par\nointerlineskip
\setlength{\abovecaptionskip}{6pt}
\caption{First-attempt native solve-time reduction relative to fixed \(w=1\) on matched joint-selected hierarchies. (a) \(60^3\) diffusion. (b) \(60^3\) variable-flow diffusion--advection. Open points show results from the six final learner states; diamonds and whiskers give their mean and sample standard deviation. Each value is computed from summed costs on 100 problems. Both panels use the same horizontal scale; all first-attempt times are included regardless of convergence.}
\label{fig:matched_policy_savings}
\end{figure}

Figure~\ref{fig:matched_policy_savings} compares the four alternatives with fixed \(w=1\). On diffusion, the global best and per-problem best fixed policies reduce mean first-attempt native solve time by \(34.0\%\) and \(35.8\%\), respectively. The period-two and frozen-LSTDQ policies give reductions of \(41.2\%\) and \(39.6\%\). On variable-flow diffusion--advection, the corresponding reductions are \(23.2\%\), \(30.9\%\), \(-3.2\%\), and \(23.6\%\), respectively.
Thus, on these matched SPD diffusion hierarchies, both varying-weight policies outperform the tested fixed-weight references on average in first-attempt native solve time. Together with the runtime decomposition in Figure~\ref{fig:runtime_breakdown}, this supports the 
benefit of the
joint strategy: \emph{low-cost hierarchy construction is paired with relaxation control that improves the subsequent solve time.}

\begin{table}[!htbp]
\centering
\caption{Frozen-LSTDQ first-attempt native solve-time reduction relative to matched-hierarchy references (\%). Reductions are mean \(\pm\) sample standard deviation across the six final learner states for each \(60^3\) problem family. Positive values favor LSTDQ. Failure rates are computed over 600 first attempts per policy and family.}
\label{tab:matched_policy}
\footnotesize
\begin{tabular}{@{}lrrrr@{}}
\toprule
& \multicolumn{2}{c}{Diffusion} & \multicolumn{2}{c}{Variable flow} \\
\cmidrule(lr){2-3}\cmidrule(l){4-5}
Policy/reference & Time reduction & Failure rate & Time reduction & Failure rate \\
\midrule
Fixed $w=1$ & $39.6\,\pm\,1.5$\% & $0.5$\% & $23.6\,\pm\,5.8$\% & $5.2$\% \\
Global best fixed & $8.5\,\pm\,2.1$\% & $0.2$\% & $0.4\,\pm\,6.1$\% & $1.0$\% \\
Per-problem best fixed & $5.9\,\pm\,2.0$\% & $0.0$\% & $-10.6\,\pm\,7.0$\% & $0.2$\% \\
Periodic $(2.85,1.10)$ & $-2.8\,\pm\,3.8$\% & $0.3$\% & $25.8\,\pm\,6.0$\% & $43.0$\% \\
Frozen LSTDQ & --- & $0.2$\% & --- & $0.3$\% \\

\bottomrule
\end{tabular}
\end{table}

On diffusion, Table~\ref{tab:matched_policy} shows that frozen LSTDQ reduces mean first-attempt native solve time by \(8.5\%\) relative to the global best fixed weight and by \(5.9\%\) relative to the per-problem best fixed weight. Frozen LSTDQ has \(2.8\%\) higher mean first-attempt native solve time than the prescribed period-two schedule. Thus, varying weights improves on the tested fixed-weight references, while a simple SPD-derived schedule captures much of this benefit without learning.

On variable-flow diffusion--advection, the period-two schedule fails to converge within the prescribed cycle limit on \(43.0\%\) of the test problems (258 of 600), predominantly at small diffusion coefficients \((\epsilon<0.005)\). In contrast to the SPD diffusion results, frozen LSTDQ achieves a \(25.8\%\) reduction in native solve time relative to the
period-two schedule, with a first-attempt failure rate of \(0.3\%\). It achieves a slight mean solve-time reduction of \(0.4\%\) relative to the global best fixed reference, but a reduction of \(-10.6\%\) relative to the per-problem best fixed reference. The failure rates of these two fixed references are \(1.0\%\) and \(0.2\%\), respectively. The robustness advantage therefore holds relative to fixed \(w=1\) and the prescribed period-two schedule; the results do not establish consistent superiority over the optimized per-problem fixed-weight references across different PDE families.

\begin{figure}[!htbp]
\centering
\includegraphics[width=\linewidth]{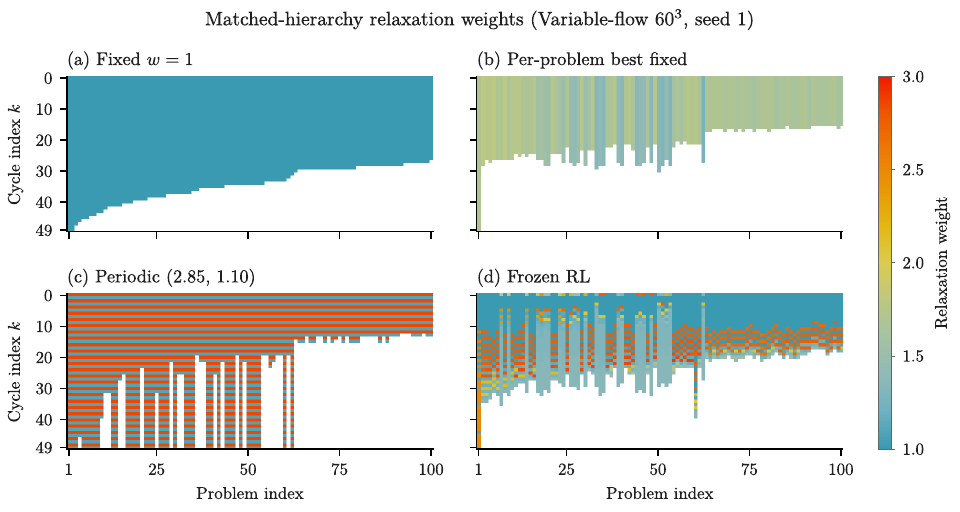}
\par\nointerlineskip
\setlength{\abovecaptionskip}{6pt}
\caption{Recorded weights for fixed \(w=1\), the per-problem best fixed weight, period-two \((2.85,1.10)\), and frozen LSTDQ on identical joint-selected hierarchies. The illustration uses seed 1. All panels contain the same 100 problems, ordered by decreasing first-attempt cycle count for \(w=1\), with case ID breaking ties. Cycle zero is at the top; white cells follow primary-attempt termination. Recovery cycles are omitted.}
\label{fig:matched_policy_weights}
\end{figure}

Figure~\ref{fig:matched_policy_weights} shows the different forms of relaxation control in one illustrative run on variable-flow diffusion--advection. Columns are ordered by the fixed-\(w=1\) first-attempt cycle count, so their left-to-right progression does not represent continued learning. The per-problem best fixed weight in (b) changes between problems but remains constant within each solve. The period-two policy in (c) follows the same prescribed alternating sequence on every problem, with limited transferability to this variable-flow family. In panel (d), frozen LSTDQ often maintains weights near \(w=1\) for several consecutive cycles before switching to alternating choices whose onset and duration vary. Among the 37 problems on which the period-two schedule fails to converge in this illustrative run, frozen LSTDQ converges on 36. On many of these nonsymmetric problems, it adopts a substantially different relaxation strategy, maintaining relatively constant low weights over extended portions of the solve rather than using strict high--low alternation. These weights are similar to the per-problem best fixed weights. For example, at displayed problem index 50, frozen LSTDQ holds \(w=1.30\) for 15 consecutive cycles, close to the per-problem best fixed weight of \(w=1.35\) for that problem. The displayed run illustrates these qualitative differences and highlights the adaptive nature of our hierarchy-aware LSTDQ policy across different problems.

\section{Conclusion}
\label{sec:conclusions}

We developed an online autotuner for \hypre's BoomerAMG that combines LinUCB-based AMG setup selection with LSTDQ-based cycle-level relaxation control, without offline pretraining. Across the tested problem sequences, jointly adapting the AMG hierarchy and relaxation weights reduces total recorded runtime relative to the default, with generally larger gains later in each sequence.
The experiments show the separate benefits of bandit-based setup tuning and RL-based solve control, as well as their coordinated adaptation during online learning. The hierarchy-aware nonstationary relaxation policy improves solve performance, while setup learning adapts hierarchy selection to the evolving solve policy. 

Overall, the results suggest that online autotuning is most useful when effective solver choices cannot be prescribed reliably in advance. Our reliance on measured AMG runtime also makes the learning process sensitive to timing variability and execution conditions. Future work will consider broader problem families, improved exploration, and large-scale parallel settings.

\section*{Acknowledgment}
The authors used GPT 6.1 to assist with manuscript drafting and editing, literature review, and theoretical development, and Codex to assist with coding. The authors take full responsibility for the accuracy and integrity of the work.

\bibliographystyle{siamplain}
\bibliography{references}

\end{document}